\documentclass{TP} 

\usepackage{amsthm,amssymb}
\usepackage[leqno]{amsmath}
\usepackage{amscd}
\usepackage[all]{xy}
\usepackage{graphicx}

\newtheorem{thm}{Theorem}
\newtheorem{lem}[thm]{Lemma}
\newtheorem{prop}[thm]{Proposition}
\newtheorem{cor}[thm]{Corollary}

\theoremstyle{definition}

\theoremstyle{remark}
\newtheorem*{rmk}{Remark}

\newenvironment{pf}{ \begin{proof} }{ \end{proof} }

\def\BbK{\mathbb{K}}
\def\Q{\mathbb{Q}}
\def\R{\mathbb{R}}
\def\Z{\mathbb{Z}}
\def\C{\mathbb{C}}

\DeclareMathOperator{\Ext}{\mathsf{Ext}}
\DeclareMathOperator{\Char}{\mathsf{char}}
\DeclareMathOperator{\tw}{\mathsf{tw}}

\DeclareMathAlphabet\EuScript{U}{eus}{m}{n}
\SetMathAlphabet\EuScript{bold}{U}{eus}{b}{n}

\DeclareFontFamily{U}{eus}{\skewchar\font'60}%
\DeclareFontShape{U}{eus}{m}{n}{<-6>eusm5<6-8>eusm7<8->eusm10}{}
\DeclareFontShape{U}{eus}{b}{n}{<-6>eusb5<6-8>eusb7<8->eusb10}{}

\begin{document}

\title{Fukaya categories of the torus and Dehn surgery}

\author{Yank{\i} Lekili}
\author{Timothy Perutz}

\begin{abstract}
This paper is a companion to the authors' forthcoming work extending  Heegaard Floer theory from closed 3-manifolds to compact 3-manifolds with two boundary components via quilted Floer cohomology. We describe the first interesting case of this theory: the invariants of 3-manifolds bounding $S^2\amalg T^2$, regarded as modules over the Fukaya category of the punctured 2-torus. We extract a short proof of exactness of the Dehn surgery triangle in Heegaard Floer homology.  We show that  $A_\infty$-structures on the graded algebra $A$ formed by the cohomology of two basic objects in the Fukaya category of the punctured 2-torus are governed by just two parameters $(m^6,m^8)$, extracted from the Hochschild cohomology of $A$. For the Fukaya category itself, $m^6\neq 0$. 
\end{abstract}

\maketitle
\section{Introduction}
This article is an offshoot of the authors' forthcoming work \cite{LP}. In that paper we will combine a detailed geometric examination of the Lagrangian correspondences between symmetric products of Riemann surfaces, studied by the second author in \cite{Per}, with the $A_\infty$ quilted Floer theory of Ma'u--Wehrheim--Woodward \cite{MWW} and the functoriality principle of \cite{WWfunct} (see also \cite{LL}). By doing so, we will extend the package of Heegaard Floer cohomology invariants \cite{OS} from closed 3-manifolds to compact 3-manifolds with boundary. To be precise, we construct invariants for compact, oriented, connected 3-manifolds with precisely two boundary components, marked as `incoming' and `outgoing'. When these are both spherical, our invariants capture the Heegaard Floer cochains of the capped-off 3-manifold. We refer to Auroux's work \cite{Aur} for the relationship of this theory to bordered Heegaard Floer theory \cite{LOT}.
 
The format of our invariants is alarmingly abstract: they take the form of $A_\infty$-functors between $A_\infty$-categories associated with the boundary surfaces, satisfying a composition law under sewing of cobordisms. Enthusiasts for extended TQFT will approve of this formulation, but geometric topologists will want to know how to extract topological information from it.  

In this article, we examine the next-to-simplest case of the theory by applying it to manifolds $Y^3$ with incoming boundary of genus $0$ (which we cap off to form $\bar{Y}$) and outgoing boundary of genus $1$. The relevant $A_\infty$-categories are certain versions of the Fukaya category of a symplectic 2-torus $T$ with a distinguished point $z$. The simplest version of the invariant for $Y$ is 

\begin{center}
{\bf an $A_\infty$-module $\widehat{\EuM}_Y$ over the Fukaya category \\$\wh{\EuF}(T_0)$ of exact, embedded curves in $T_0:=T\setminus \{z\}$.}
\end{center}
This module evaluates on each object $X$ (which is a circle $X\subset T$ with an exactness constraint and certain decorations) to give a cochain complex $\wh{\EuM}_Y(X)$. This complex is  quasi-isomorphic to the Heegaard Floer cochains $\wh{CF^*}(\bar{Y}\cup_T U_X)$, where $U_X$ is the solid torus in which the circle $X$ bounds a disc. The different objects $X$ correspond to different Dehn fillings of $Y$.

We illustrate our theory by showing how it leads to a proof of the Dehn surgery exact triangle in Heegaard Floer cohomology  \cite{OSapp} (Theorem \ref{exacttri}, Corollary \ref{surgery}). By working with the modules $\widehat{\EuM}_Y$ rather than the Heegaard Floer cochain complexes, one can legitimately work on the genus 1 boundary rather than on the $g$-fold symmetric product of a genus $g$ Heegaard surface. This makes the proof technically straightforward and also makes the signs, on which the proof depends, transparent.

To better understand the nature of our invariants, we need to understand the structure of $\wh{\EuF}(T_0)$. This category is closely related to the one studied in $\cite{PZ}$ (see \cite{AS} for a deformation-theoretic approach to the latter category), but the differences are significant. It follows from the surgery exact triangle (Theorem \ref{exacttri}) that $\wh{\EuF}(T_0)$ (whose objects we declare to be exact, oriented Lagrangians with \emph{non-trivial} spin structures) is generated by two objects---curves $a$ and $b$ that meet  transversely at a point, generating a full $A_\infty$-subcategory $\EuA$. That is, the inclusion $\EuA\subset \wh{\EuF}(T_0)$ induces a quasi-equivalence of the triangulated envelopes $\tw \EuA \to \tw \wh{\EuF}(T_0)$ (`$\tw$' stands for twisted complexes; see \cite[(3l)]{SeiBook}). This in turn induces an equivalence of triangulated categories between the derived categories $D\EuA=H^0(\tw\EuA)$ and $D\wh{\EuF}(T_0)=H^0(\tw \wh{\EuF}(T_0)$).  

The $A_\infty$-structure of $\EuA$ can be transferred to an $A_\infty$-structure on the cohomology algebra $A=H^*\EuA$. This structure is governed by the Hochschild cochain complex $CC^*(A,A)$, with its bigrading and Gerstenhaber bracket. We show that---over a field $\BbK$ in which 6 is invertible---gauge-equivalence classes of $A_\infty$-structures in $A$ determine and are determined by two parameters, $m^6\in
HH^2(A,A)^{2-6}\cong \BbK$ and $m^8\in HH^2(A,A)^{2-8}\cong \BbK$. 
The `moduli space' $HH^2(A,A)^{-4}\times HH^2(A,A)^{-6}$ is parametrized by the invariants $(m^6,m^8)$ of dg algebras associated with a pair of sheaves---the structure sheaf and the skyscraper at $[0:1:0]$---on the Weierstrass curves $y^2z=4x^3-pxz^2-qz^3$ (Prop. \ref{Weierstrass}) for $(p,q)\in \BbK^2$.
The value of $m^6(\EuA)$ is non-zero (Theorem \ref{m6}). Our expectation is that $\EuA$ is quasi-isomorphic to the dg algebra for a \emph{nodal} cubic.\footnote{We have confirmed this expectation, but will write up our argument elsewhere.}

\paragraph{Note.} The topics we have chosen do not rely on detailed information about our general theory. Most of the theorems proved in this article are independent of it; the sole exception, Corollary \ref{surgery}, invokes a general feature of the theory. The section on Hochschild cohomology can be read on its own. 

We invoke the definition of the Fukaya category from \cite{SeiBook}. We refer to that book for points of homological algebra, but recommend \cite[Section 3]{SeiQuartic} as a briefer alternative reference.

\section{Fukaya categories} A \emph{Liouville domain} is a compact, even-dimensional manifold $M^{2n}$ with boundary, equipped with a 1-form $\theta$ such that $\omega:=d\theta$ is everywhere non-degenerate, and such that the Liouville vector field $\lambda$, characterized by the equation $\iota(\lambda)\omega=\theta$, points outwards along the boundary. An \emph{exact Lagrangian submanifold} is a closed, embedded $n$-submanifold $L\subset \interior(M)$ such that $\theta|_L$ is an exact 1-form.

The \emph{Fukaya category} $\EuF(M)$ of a Liouville domain $M$ with $c_1(TM)=0$, as constructed in \cite{SeiBook}, is an $A_\infty$-category, linear over a field $\BbK$. Its objects are arbitrary exact Lagrangian submanifolds equipped with spin-structures and gradings. For objects $X_0$ and $X_1$ of $\EuF(M)$,  the morphism-space $\hom(X_0,X_1)$ is a Floer cochain complex $CF^*(\phi(X_0),X_1)=\BbK^{\phi(X_0)\cap x_1}$, where $\phi$ is the time-1 diffeomorphism of the Hamiltonian vector field for a function $H=H_{X_0,X_1}$, drawn from a look-up table. The first structure map $\mu^1\colon \hom(X_0,X_1)\to \hom(X_0,X_1)$ is Floer's coboundary map. The structure constants of the higher composition maps $\mu^d$ defining the $A_\infty$-structure are counts of inhomogeneous pseudo-holomorphic polygons. Up to quasi-isomorphism, the whole structure is an invariant of $(M,[\theta])$, where $[\theta]$ is the class of $\theta$ modulo closed 1-forms supported in $\interior(M)$.

The case of surfaces $M^2$ is special because one can use the Riemann mapping theorem to describe moduli spaces of holomorphic polygons, and thereby describe some of the structure maps in combinatorial terms.  We say `some' because when $\mu^d$ is applied to a sequence of curves that includes repetitions, the na\"ive moduli spaces tend to be irregular, and perturbations are required. We next give a partial combinatorial formulation, following \cite[section 13]{SeiBook} and \cite{SeiHMS2}.

\paragraph{The 1-pointed 2-torus.} Let $T$ be a closed Riemann surface of genus 1 equipped with a basepoint $z$---a complex elliptic curve. Choose a hermitian metric in the holomorphic line bundle $\EuL:=\mathcal{O}_T(z)$, and a unitary connection $A\in  \Omega^1(S\EuL, i\R)$ such that the curvature form $\omega:= (i/2\pi) F_A$ is an area form, positive with respect to the complex orientation.  Choose also a tangent line field $\lambda$. Write $T_0$ for the Riemann surface $T\setminus \{z\}$. The restriction $\EuL|_{T_0}$ is tautologically trivialized; in this trivialization, we can write $A=d- 2\pi i \theta$, where $d\theta=\omega$. Notice that $\theta$ restricts to a small loop encircling $z$ as a contact 1-form; the Liouville field points towards $z$ along this loop.

Define a \emph{marked exact Lagrangian} to be a pair $(\gamma,\star)$, where $\gamma$ is an oriented circle embedded in $T_0$ as an exact Lagrangian submanifold (that is, $\int_\gamma\theta=0$), and $\star\in \gamma$.  A \emph{grading} on a marked exact Lagrangian is a homotopy from $\lambda|_\gamma$ to $T\gamma$ inside the space of line fields $\mathbb{P}(TT_0|_\gamma)$. Gradings which are homotopic rel endpoints are considered identical, which has the consequence that the gradings for $\gamma$ form a $\Z$-torsor.

Non-contractible embedded curves in $T_0$ are determined, up to homotopy, by their slope, which is an arbitrary element of $\Q\cup \{\infty\}$. In each homotopy class there is an exact representative (take an initial representative $\gamma_0$ and move it by an isotopy of flux $-\int_{\gamma_0}\theta$), unique up to isotopy through exact Lagrangians. 

\paragraph{Objects.} The objects of our Fukaya category $\wh{\EuF}(T_0)$ will be marked, exact graded Lagrangians. The `hat' is there firstly because this category is relevant to the `hat' version of Heegaard Floer cohomology, and secondly because $\wh{\EuF}(T_0)$ is not quite the same as the usual exact Fukaya category $\EuF(T_0)$, but is a full subcategory of it.  The distinction is the following.  The basepoint $\star$ determines a spin-structure on $\gamma$.  A spin structure amounts to a double covering $\tilde{\gamma} \to \gamma$, and we declare this double covering to be trivial over $\gamma \setminus \{\star\}$, and to interchange the sheets at $\star$. \emph{We do not allow $\gamma$ to carry the trivial spin structure.}

\paragraph{Morphisms.} Fix the base field $\BbK$. Take a pair of objects $(X,Y)$, where $X=(\gamma_X,\star_X)$ and $Y=(\gamma_Y,\star_Y)$, and pick a Hamiltonian function $H_{X,Y}$ whose flow $\phi^t_{X,Y}$ generates a family of curves $\{ \gamma_X^t\}_{t\in [-1,1]}$, such that $\gamma_X^t \pitchfork \gamma_Y$ for all $t\in [-1,0)\cup(0,1]$. Because transversality is maintained, the combinatorial pattern of intersection points is unchanged for $t\in (0,1]$. The points $\star_X$ and $\star_Y$ should not occur as intersection points $\gamma_X^t \cap \gamma_Y$ for any $t\in (0,1]$. We define $\hom(X,Y)$ to be the Floer cochain complex $CF^*(\gamma_X^1,\gamma_Y)$.

In more detail, each point $x\in \gamma_X^1 \cap \gamma_Y$ has a Maslov index $i(x)\in \Z$ defined using the gradings \cite[section 13]{SeiBook}. The Floer complex is
\begin{equation} 
CF^*(\gamma_X^1,\gamma_Y) = \bigoplus_{x\in \gamma_X^1\cap\gamma_Y} {\BbK x}, \quad
\deg(x) = i(x). 
\end{equation}
The differential
$d\colon CF^*(\gamma_X^1,\gamma_Y) \to CF^{*+1}(\gamma_X^1,\gamma_Y)$ is  defined by linearity and
\begin{equation}
d x = \sum_{y\in \gamma_X^1\cap \gamma_Y}\sum_{u\in \EuM(x,y)}{ \sigma(u)y}. 
\end{equation}
Here $\EuM(x,y)$ is the set of homotopy-classes of \emph{immersed bigons} in $T_0$. By an immersed bigon, we mean a homotopy class of smooth, orientation-preserving immersions $u\colon \bar{D^2}\setminus \{ 1,-1 \} \to T_0$ such that $u(\{ e^{i\theta}: \theta \in (0,\pi) \}) \subset \gamma_Y$ and $u(\{ e^{i\theta} :\theta\in (\pi, 2\pi) \})\subset \gamma_X^1$, where $u$ extends smoothly to a map $\bar{u}\colon \bar{D^2}\to T_0$ with $\bar{u}(1)=x$ and $\bar{u}(-1)=y$. Furthermore, the image of the map should have \emph{convex corners} at $x$ and $y$. The set $\EuM(x,y)$ is finite, and to each $u\in \EuM(x,y)$ one can assign a sign $\sigma(u)\in \{ \pm 1\}$. This is given by $\sigma(u) = (-1)^{r+s}$, where $r$ is $0$ if the upper boundary of $u$ traverses $\gamma_Y$ in the positively oriented direction, and $1$ otherwise; and $s$ is the number of stars ($\star_X$ or $\star_Y$) encountered by $u$ on its boundary.

One has $d^2=0$, so  $CF^*(\gamma_X^1,\gamma_Y)$ is a cochain complex, and one can form the cohomology $\Hom(X,Y)=HF^*(\gamma_X^1,\gamma_Y):=H(CF^*(\gamma_X^1,\gamma_Y))$. When $\gamma_X$ and $\gamma_Y$ are exact-isotopic,  $\Hom(X,Y)\cong H^*(\gamma_X;\BbK)$. 

\paragraph{Structure maps.}
The $A_\infty$ structure maps of $\wh{\EuF}(T_0)$,
\begin{equation}\label{ainf} 
\mu^d \colon \hom(X_{d-1},X_d) \otimes \dots \otimes  \hom(X_0,X_1) \to \hom(X_0,X_d)[2-d] \end{equation}
are in general defined via solutions to an inhomogeneous Cauchy--Riemann equation, but in certain cases one can obtain consistent and regular moduli spaces of genuine holomorphic maps into $T_0$. We summarize what we need, referring to \cite[Chapter 2]{SeiBook} and \cite{SeiHMS2} for further details.

Write $X_i=(\gamma_i,\star_i)$, and suppose that all the $X_i$ are drawn from a fixed, finite set $O$ of objects such that the curves underlying two distinct objects intersect transversely (and not at the $\star$-points) and all triple intersections are empty. In that case, we can suppose that the functions $H_{X_i,X_j}$ are zero when $X_i\neq X_j$.

The most pleasant instance of (\ref{ainf}) occurs when any two of the objects $X_i$ in the sequence are distinct elements of $O$ (hence transverse as curves).  In that case, consider intersection points $y_i\in \gamma_{i-1}\cap \gamma_i$ for $i=1,\dots,d$. One has 
\[ \mu^d(y_d,\dots, y_1)=\sum_{y_0\in  \gamma_d\cap \gamma_0} {n(y_0;y_1,\dots,y_d) y_0},\] 
where  $n(y_0;y_1,\dots,y_d)$ is a signed count of immersed polygons. 
To be precise, the relevant polygons are smooth  immersions of $\bar{D}\setminus \{ e^{2\pi i k/(d+1)}: k\in \Z\}$ into $T_0$. The immersion must preserve orientation, must extend continuously to a map on $\bar{D}$ sending $e^{2\pi i k/(d+1)}$ to $y_k$, and must map the boundary interval $\{ e^{2\pi i t/(d+1)}: t\in (k-1,k)\}$ to $\gamma_k$. It must also have convex corners. What we count are the homotopy classes of such immersions. 

The sign attached to an immersion is $(-1)^{q+r+s}$. As before, $s$ is the number of stars on the boundary; $q$ is $i(y_0)+i(y_d)$ if the polygon travels along $\gamma_d$ in the negative direction with respect to the orientation, and $0$ otherwise; and $r$ is the sum of degrees $i(y_k)$ over those $k\in \{1,\dots, d-1\}$ (not $k=0$) such that the polygon travels along $\gamma_k$ in the negative direction.

Now suppose that the sequence $(\gamma_0,\gamma_1,\dots,\gamma_d)$ includes precisely one repetition, occurring between cyclically adjacent curves (i.e., $X_k=X_{k+1}$ for some $k \in \Z/d$) with all other pairs distinct, hence transverse. In this case, we replace $\gamma_k$ by $\gamma_k'$, the image of $\gamma_k$ under the time $(-1)$ Hamiltonian diffeomorphism $\phi^{-1}_{H}$ for the function $H=H_{X_k,X_k}$ associated with the pair $(X_k,X_k)$. One then applies the above recipe to the transverse sequence of objects with $X_k$ replaced by $X_k'$. 

As stated this does not quite make sense because after this replacement, the intersection points are between the wrong curves.  However, when $X_k\neq X_j$, one has a canonical isomorphism of cochain complexes $CF(\gamma_k', \gamma_j)\cong \hom(X_k, X_j)$ arising from the bijection
$\gamma_k' \cap \gamma_j \cong \gamma_k \cap \gamma_j$ coming from the flow $\phi^t_H$. One also has a canonical isomorphism $CF(\gamma_k,\gamma_k')\cong CF(\phi^1_H(\gamma_k),\gamma_k))$, induced by 
$\phi^1_H$. One uses these isomorphisms to make sense of this formulation.

\paragraph{Versions including the basepoint.}

As a variant on the construction of $\wh{\EuF}(T_0)$, one can construct the
$A_\infty$-category $\EuF_\infty(T,z)$ with the same objects, and morphism
spaces $\hom_\infty(X,Y)=\hom_{\EuF(T_0)}(X,Y)\otimes \BbK[U^{-1},U]]$, where
$U$ has degree 0. \footnote{In Heegaard Floer theory, one usually declares $U$
to have degree 2. The genus 1 case is anomalous. The stabilization isomorphism
relating Heegaard complexes computed in genera 1 and 2 is not
degree-preserving, but it does respect a natural ``geometric grading'', refining the $\spinc$ grading, by the
$\Z$-set of homotopy classes of oriented 2-plane fields on the 3-manifold.} The structure maps are constructed
in just the same way as before, except that the immersed polygons are now
allowed to pass through $z$, and $n(y_0;y_1,\dots,y_d)$ counts such a polygon
$u$ with a weight $U^m$, where $m$ is the multiplicity of $z$ in $y$. (This
procedure is familiar from Heegaard Floer theory). Since the multiplicities are
non-negative, one can construct $A_\infty$-categories $\EuF_+(T,z)$ and
$\EuF_-(T,z)$ whose morphism spaces are \begin{align*}
	&\hom_+(X,Y)=\hom_{\EuF(T_0)}(X,Y)\otimes \BbK[[U]],\\
	&\hom_-(X,Y)=\hom_{\EuF(T_0)}(X,Y)\otimes \BbK[U^{-1},U]]/\BbK[[U]].
\end{align*} $\wh{\EuF}(T_0)$ and the $\pm$-categories are all three determined
by the $\infty$ version, together with its filtration by the subcategories with
hom-spaces $\hom_{\EuF(T_0)}(X,Y)\otimes U^k\BbK[[U]]$.  \section{Generalized
Lagrangians and Fukaya-modules} There is a contravariant Yoneda embedding of
$\EuF(M)$ into the $A_\infty$-category $Mod\, \EuF(M)$ of right
$\EuF(M)$-modules, that is, $A_\infty$ functors from $\EuF(M)$ to the dg
category of cochain complexes. The embedding maps an object $X$ of $\EuF(M)$ to
the module $\EuY_X$ which assigns to each object $X'$ the cochain complex
$\EuY_X(X'):=\hom(X',X)$ \cite{SeiBook}.

There is also an intermediate category $\EuF^\sharp(M)$, the \emph{extended} Fukaya category, whose objects are finite sequences of Lagrangian correspondences (of a constrained kind) between Liouville domains, where the sequence begins at $\{pt.\}$ and ends at $M$. The morphism-spaces are quilted Floer cochain complexes, as in \cite{WWquilt}. Making this precise is the substantial task of \cite{MWW}.  One has embeddings
\begin{equation} 
 \EuF(M) \to \EuF^\sharp(M) \stackrel{\EuY^\sharp}\to Mod\,\EuF(M)
 \end{equation}
factoring $\EuY$; the first functor is an obvious inclusion, $\EuY^\sharp$ another Yoneda embedding. 

The construction in \cite{LP} will attach to a 3-manifold $Y$ bounding $T$ an object $L_Y$ in $\EuF^\sharp(T_0)$. It depends on additional choices, but the resulting module $\wh{\EuM}(Y ):= \EuY^\sharp_{L_Y} $ over $\wh{\EuF}(M)$ does not (up to isomorphism). The object $L_Y$ also defines a filtered module $\EuM_\infty(Y)$ over $\EuF_\infty(T,z)$, hence modules $\EuM_+(M)$ and $\EuM_-(M)$ over $\EuF_+(T,z)$ and $\EuF_-(T,z)$.

In this article, we shall use only one property of these modules, to be proved in \cite{LP}, which is as follows. Let $X=(\gamma_X,\star_X)$ be an object in $\EuF_\infty(T_0)$, and $Y(\gamma)=Y\cup_{\gamma\sim \partial D} (S^1\times D^2)$ the closed 3-manifold obtained from $Y$ by Dehn-filling $\gamma$. 

($\dagger$) \emph{The $\BbK[U^{-1},U]]$-linear cochain complex $(\EuM_\infty Y)(X)$ is quasi-isomorphic to the completion at $U$ of the Heegaard Floer cochains $CF^*_\infty(Y(\gamma))$ defined in \cite{OS}, compatibly with the natural filtrations of these complexes.}  

\section{The surgery exact triangle} 
In an $A_\infty$-category $\EuC$, one has a notion of an \emph{exact triangle} $X\stackrel{a}{\to} Y \stackrel{b}{\to} Z \stackrel{c}{\to} X[1]$: a trio of morphisms (only their classes $[a]$, $[b]$ and $[c]$ in $H^0\EuC$ matter) which is transformed by the Yoneda embedding $\EuY \colon \EuC \to mod\,\EuC$ into a triangle of $\EuC$-modules $\EuY_X \to \EuY_Y \to \EuY_Z \to \EuY_X[1]$, isomorphic in $H (mod\,\EuC)$ to the standard triangle $\EuY(X)\to \EuY(Y) \to \cone(a) \to \EuY(X)[1]$ 
associated with the mapping cone of $a$ \cite[(3f)]{SeiBook}.
\\
\begin{figure}\label{ExactTriangle}
\begin{center}
\includegraphics[scale=0.67]{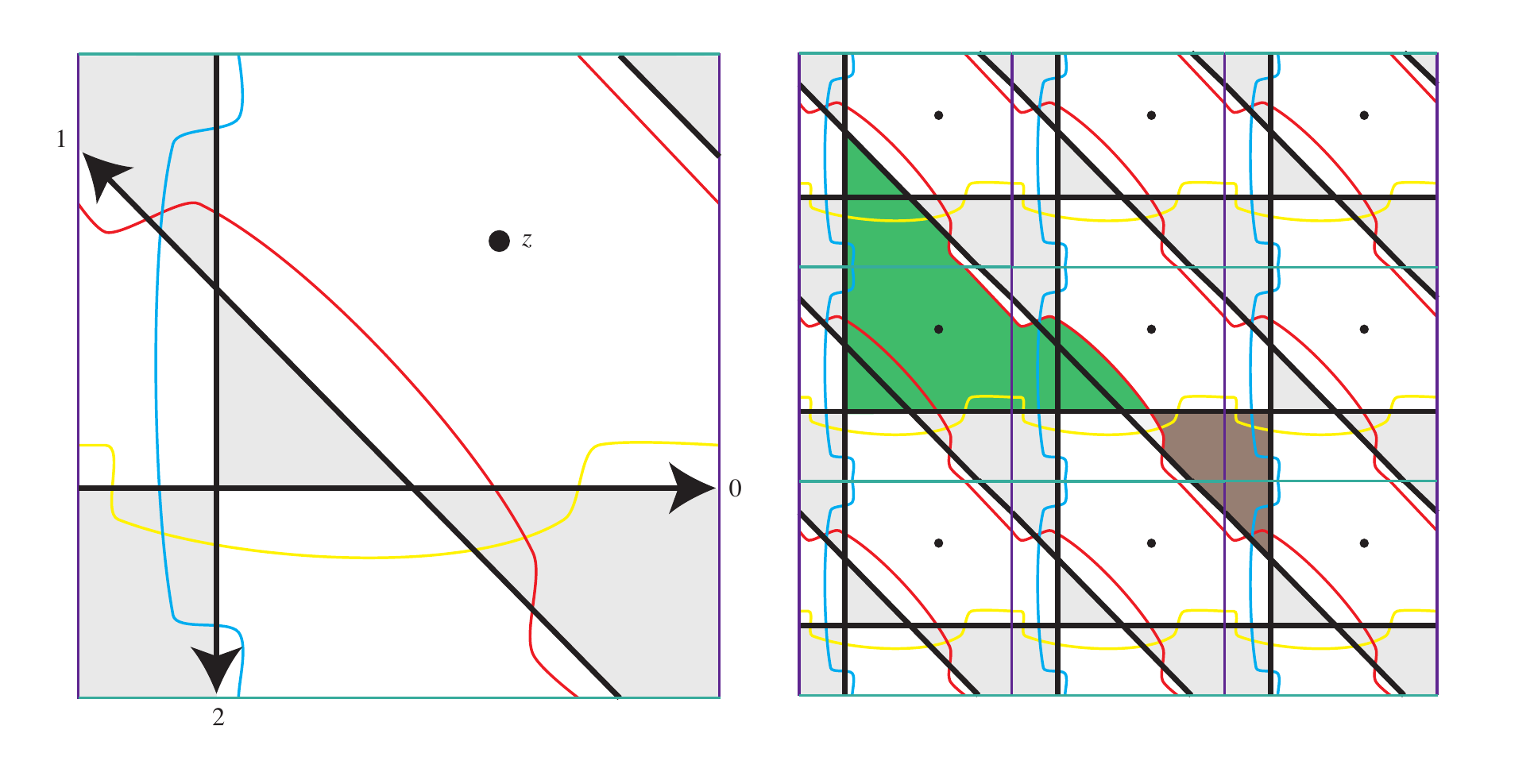}
\caption{Left: the oriented curves $\gamma_0$, $\gamma_1$ and $\gamma_2$ in the torus, with the basepoint $z$ placed in the hexagonal part of their complement. Hamiltonian deformations of these curves, used to define endomorphisms in the Fukaya category, are shown in color. Right: Immersed polygons in the torus viewed as  embedded polygons in the universal cover.}
\end{center}
\end{figure}

\begin{thm}\label{exacttri}
Let $X_0=(\gamma_0,\star_0)$, $X_1=(\gamma_1,\star_1)$ and 
$X_2=(\gamma_2,\star_2)$ be three objects in $\wh{\EuF}(T_0)$. Suppose that any two of the underlying curves, $\gamma_i$ and $\gamma_j$, say, intersect at a single point $e_{ij}$, transversely, and that the curves are oriented so that
\[ (\gamma_0\cdot \gamma_1) = (\gamma_1\cdot \gamma_2) =(\gamma_2\cdot \gamma_0) = 1,
\]
as pictured in Figure 1 (left).  Then one has an exact triangle in $\EuF_\infty(T,z)$, 
\begin{equation} \label{ex tri}
X_0 \stackrel{-ue_{01}}{\to} X_1 \stackrel{-u e_{12}}{\to} X_2
\stackrel{ -u e_{20}}{\to} X_0[1], 
\end{equation}
where $u$ is Euler's generating function for partitions,
\[ u= \prod_{m> 0}{(1-U^m)^{-1}}\in \BbK[[U]]^\times.\] 
This triangle remains exact in $\EuF_+(T,z)$ and in $\EuF_-(T,z)$, and specializes (by putting $U=0$ in the $+$ version) to an exact triangle in $\wh{\EuF}(T_0)$.
\end{thm}

The exactness of (\ref{ex tri}) in $\wh{\EuF}(T_0)$ is a special case of Seidel's exact triangle \cite[Section 17j]{SeiBook} (see also \cite{SeiLES}). The full statement is new, the $u$-factors in particular, but cf. \cite{PolMassey}.

\begin{pf}
The curves depicted as bold black lines in Figure 1 (left) represent $\gamma_0$, $\gamma_1$ and $\gamma_2$; the orientations are shown as arrowheads. The gradings will turn out to be immaterial. In $\EuF_\infty(T,z)$, one has
\[ \hom(X_i,X_j) = \mathbb{K}[U^{-1},U]] \cdot e_{ij}, \quad i\neq j.  \]

We choose Floer data on $X_i$ such that $\hom(X_i,X_i)=CF(X_i',X_i)$, where $X_i'$ is the object whose underlying curve $\gamma_i'$ is the small Hamiltonian push-off of $\gamma_i$ depicted in Figure 1 (left).  Let $\EuT$ be the full subcategory of $\EuF_\infty(T,z)$ with objects $X_0$, $X_1$ and $X_2$; notice that, with our choice of Floer data, $\mu^1_{\EuT}=0$.  One can recognize exact triangles using the criterion of \cite[Lemma 3.7]{SeiBook}, which is particularly straightforward since $\mu^1_\EuT=0$. To prove that (\ref{ex tri}) is exact in $\EuT$ (or equally in $\EuF_\infty(T,z)$), it is sufficient that
\begin{equation}\label{criterion}
  \mu^2_\EuT(e_{01} ,e_{20}) =0 =  \mu^2_\EuT(e_{20},e_{12} ), \quad   -u^3 \mu^3_\EuT(e_{01}, e_{20},e_{12}) = e_{X_1},
\end{equation}
where $e_{X_1}$ is the unique element in $\hom(X_1,X_1)$ representing the identity in $\Hom(X_1,X_1)$; and that for each object $Z$ of $\EuF_\infty(T,z)$, the cochain complex
\begin{equation}\label{sum of three} 
\hom(Z, X_2)[1] \oplus \hom(Z, X_0)[1] \oplus \hom(Z,X_1)  
\end{equation}
is acyclic with respect to the differential
\[ d = \left[   \begin{array}{ccc}
    \mu^1 & 0 & 0 \\ 
    -\mu^2(ue_{20},\cdot) & \mu^1 & 0 \\ 
    \mu^3(ue_{01}, ue_{20},\cdot) & -\mu^2(ue_{01},\cdot)  & \mu^1 \\ 
\end{array}\right]. \]

We compute $\mu^2_\EuT$ by examining immersed triangles with convex corners. For the composite $\mu^2_\EuT(e_{01} ,e_{20})$, there are two embedded triangles, shaded in the left-hand figure,
neither of them intersecting $z$, contributing a leading-order term $(\epsilon_{1,0} +\epsilon_{2,0})e_{21}$ for certain signs 
$\epsilon_{1,0},\epsilon_{2,0} \in \{ \pm 1\}$. There are further immersed triangles, whose lifts to the universal cover of $T$ are visible in Figure 1 (right). For each integer $p \geq 0$, there are exactly
two immersed triangles whose sides wrap around $\gamma_0$ more than $p$ times but less than $p+1$ times. These triangles pass through $z$ at $p(p+1)/2$ points.  Hence
\begin{equation} \mu^2_\EuT(e_{01} , e_{20}) = \left( \sum_{p\geq
0}{(\epsilon_{1,p} + \epsilon_{2,p})U^{p(p+1)/2} } \right) e_{21}
\end{equation}
for signs $\epsilon_{j,p} = \pm 1$.
We can also see that
\begin{equation}\label{Massey}
\mu^3_\EuT(e_{01}, e_{20}, e_{12}) =  \left( \sum_{p\geq 0}{( (p+1)\epsilon_{3,p} +p\epsilon_{4,p}) U^{p(p+1)/2}} \right) e_{X_1} 
\end{equation}
for further signs $\epsilon_{3,p},\epsilon_{4,p}  \in \{ \pm 1\}$. The leading term $\epsilon_{3,0}$ (resp. $\epsilon_{4,1}$) comes from the embedded quadrilateral shown in brown (resp. green) on the right of the figure. The higher  $\epsilon_{3,p}$ (resp. $\epsilon_{4,p}$) terms account for larger, immersed quadrilaterals, approximately similar to the two depicted.

\begin{figure}
\label{Signs}
\begin{center}
\includegraphics[scale=0.8]{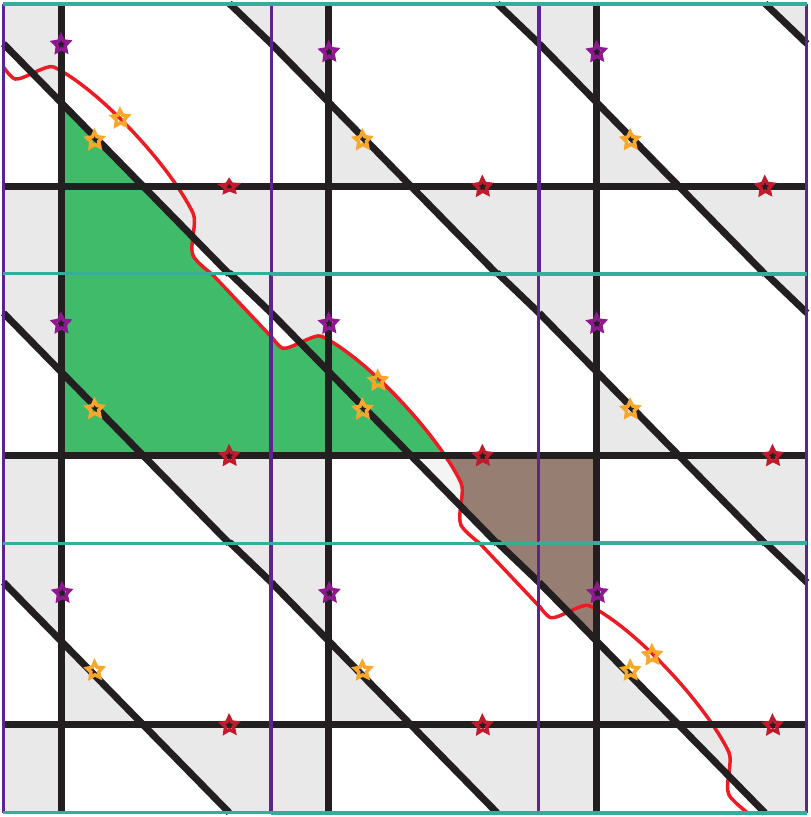}
\caption{The number of stars on the boundary of a polygon contributes to its sign.}
\end{center}
\end{figure}

Now we check the signs, using the recipe described earlier involving orientations and stars. In Figure 2 (which depicts the universal
cover), the points $\star_i$ are shown as colored stars. We read off  from Figure 2 that
\begin{align*} 
&\epsilon_{1,p} = \rho_1 (-1)^{-1+ 3p} ,&&
\epsilon_{2,p} = \rho_2 (-1)^{2+3p} ,\\
&\epsilon_{3,p} = \rho_3 (-1)^{2+3p}  && 
\epsilon_{4,p} = \rho_4 (-1)^{4+3(p-1)} 
\end{align*}
where $\rho_j$ is the sign coming from orientation-mismatches, and the power of $(-1)$ is the `star-sign' $(-1)^s$. The  $\epsilon_{1,p}$- and $\epsilon_{4,p}$-polygons have boundary orientations consistent with the given orientations of the $\gamma_i$, so $\rho_1=\rho_4=1$. The $\epsilon_{2,p}$ and $\epsilon_{3,p}$-triangles have boundary orientations opposite to the given orientations. For a convex-cornered $(d+1)$-gon with corners $y_0$ (outgoing) and $(y_1,\dots,y_d)$ (incoming), in their natural boundary order, one has the relation $i(y_0)=i(y_1)+\dots + i(y_d)+2-d$; from this we see that $\rho_2=1$ and $\rho_3=-1$. Consequently, $\mu^3_\EuT(e_{01}, e_{20},e_{12}) = -v e_{X_1}$, where
\begin{equation} v= \sum_{p\geq 0}{(-1)^p(2p+1)U^{p(p+1)/2}}. \end{equation}
One has $v= u^{-3}$: this specialization of Jacobi's triple product identity can be found in his 1829 treatise \emph{Fundamenta novae theoriae ellipticarum} (section 66).

We turn now to the verification that the complex (\ref{sum of three}) has no cohomology.  This is well-known (cf. \cite{SeiLES, OSapp}), so we shall be brief. Given $Z$, move $\gamma_1$ by a Hamiltonian isotopy which does not change the combinatorics of the pattern of intersection points between $\gamma_1$ and $\gamma_0$, nor between $\gamma_1$ and $\gamma_2$, to a curve $\gamma_1'$ which lies very close to $\gamma_0\cup \gamma_2$. The Hamiltonian isotopy produces an isomorphism $\gamma_1\to \gamma_1'$ in the Fukaya category, and the cohomology of (\ref{sum of three}) is unaffected by replacing $\gamma_1$ by $\gamma_1'$ via this isomorphism.

There is a bijection
\begin{equation}  \beta\colon (\gamma_0 \cap Z) \amalg (\gamma_2\cap Z) \cong  \gamma_1' \cap Z  \end{equation}
which sends an intersection point of $\gamma_0\cap Z$ or of $\gamma_2\cap Z$ to a nearby point in $\gamma_1' \cap Z$. Work, for now, in $\wh{\EuF}(T_0)$. There are finitely many holomorphic polygons in $T_0$ contributing to the matrix entries of the differential in (\ref{sum of three}), and up to deformation, these curves do not change as we let $\gamma_1'$ degenerate towards $\gamma_0\cup \gamma_2$. Some of these curves are `small' in the sense that their area goes to zero in this degeneration; the others are `large'.  There is a unique small triangle contributing to  $\mu^2( e_{01}, x)$, and its contribution is $\pm \beta(x)$. For $y\in \gamma_2\cap Z$, there are no small-area triangles contributing to $\mu^2( e_{20}, x)$, but there is a unique small-area quadrilateral contributing to $\mu^3(e_{01}, e_{20},y)$: its contribution is $\pm \beta(y)$.  There are no small-area bigons contributing to the diagonal entries of $d$. The homology of the `small area' part of the differential is zero, and by a filtration argument based on the symplectic action (as in \cite{SeiLES}) it follows that $H^*(d)=0$. 

Another filtration argument, this time based on the filtration by powers of $U$, shows that $H^*(d)=0$ also in $\EuF_\infty(T,z)$ and the other variants.
\end{pf}

\begin{cor}[Surgery triangle]\label{surgery} Let $Y$ be a 3-manifold bounding the torus $T$. Let $\gamma_0$, $\gamma_1$ and $\gamma_2$ be simple closed curves in $T$ so that any two of them intersect
transversely at a simple point, and oriented so that 
\[(\gamma_0 \cdot \gamma_1) =(\gamma_1 \cdot \gamma_2) = (\gamma_2\cdot \gamma_0 ) = 1.\] 
Let $Y_i$ be the
3-manifold obtained by Dehn-filling  $\gamma_i\subset T$. Then there is an exact triangle of $\EuF_\infty(T,z)$-modules, 
\[ \EuM_\infty(Y_0) \to \EuM_\infty(Y_1)\to \EuM_\infty(Y_2) \to
\EuM_\infty(Y_0)[1], \] 
and hence by the principle $(\dagger)$ an exact triangle of completed Heegaard Floer cohomology groups $HF^\bullet_\infty$. The same holds true in the $+$, $-$ and `hat'  versions of $\EuF$, $\EuM$ and $HF^*$.
\end{cor} 
Define $C$ to be the smallest set of curves on $T_0$, containing $X_0$ and $X_1$, such that given any two curves $Y_0, Y_1\in C$, if $Y_2$ is a third curve such that $(Y_0,Y_1,Y_2)$ form a `surgery triangle' as in Theorem \ref{exacttri} then $Y_2\in C$. It is easy to check that $C$ contains curves of arbitrary rational slope. Hence:
\\
\begin{cor} 
The objects $a$ and $b$ generate the triangulated $A_\infty$-category of twisted complexes $\tw  \EuF_\infty(T_0)$. They likewise generate $\tw \EuF_+(T,z)$, $\tw \EuF_-(T,z)$ and $\tw \wh{\EuF}(T_0)$.
\end{cor}
\section{Computing Hochschild cohomology}\label{HH}
\subsection{A two-object subcategory} Let $\EuA$ denote the full subcategory of $\wh{\EuF}(T_0)$ with objects $a$ and $b$. We have just seen that the induced inclusion $\tw \EuA \subset \tw \wh{\EuF}(T_0)$  is a quasi-equivalence. Let $A=H^0\EuA$ be the cohomology category of $\EuA$, a graded $\BbK$-linear category. We sometimes regard $A$ (more precisely, the direct sum of its hom-spaces) as an associative, unital $\BbK$-algebra with two distinguished idempotent elements whose sum is 1. 

The normalized Hochschild cochain complex of the unital algebra $A$, $(CC^*(A,A),\delta)$, is the direct sum of subspaces $CC^{r+s}(A,A)^s$, the $\BbK$-vector space of degree-preserving linear maps $A^{\otimes r}\to A[s]$ which vanish on monomials one of whose factors is $1$. Besides being a cochain complex with an internal grading $s$ preserved by the differential $\delta$, the Hochschild complex has its Gerstenhaber bracket which makes $(CC^{*+1}(A,A),\delta)$ a dg Lie algebra. In a precise sense, cochains $\mu^\bullet \in CC^2(A,A)$ which satisfy the Maurer--Cartan equation $\delta \mu^\bullet + \frac{1}{2}[\mu^\bullet,\mu^\bullet]=0$ parametrize $A_\infty$-structures with cohomology algebra $A$ (when 2 is invertible) \cite[(3a,e)]{SeiQuartic}. It is therefore of interest to compute the Hochschild cohomology $HH^*(A,A)=H^*(CC^*(A,A),\delta)$. This we accomplish using the interpretation of Hochschild cohomology as a right-derived functor for graded bimodule homomorphisms.

\paragraph{Structure maps in $\EuA$.} An $A_\infty$-category is \emph{minimal} if $\mu^1=0$.  In  $\EuA$, the self-homs $\hom_\EuA(a,a)$ and $\hom_\EuA(b,b)$ are defined as Floer complexes $CF^*(a^1,a)$ and $CF^*(b^1,b)$, where $a^1$ is the image of $a$ under a chosen Hamiltonian diffeomorphism, and $b^1$ is the image of $b$ under a (possibly different) Hamiltonian diffeomorphism. One can arrange that $a^1\cap a$ and $b^1 \cap b$ are transverse intersections, each consisting of two points.  By making such choices, we ensure that $\EuA$ is minimal.

There is a unique point $x\in  a\cap b$. This point represents a class $u_1\in \hom^1(a,b)$ which generates $\hom^*(a,b)$. By assigning $u_1$ degree 1, we are in effect pinning down the gradings of $a$ and $b$ up to a common shift of both $a$ and $b$. The point $x$ also represents a class $v_0\in \hom^0(b,a)$, generating $\hom^*(b,a)$. That $v_0$ has degree $0$ is an instance of Floer-theoretic Poincar\'e duality: $\deg(v_0)=1-\deg(u_1)$.  One has
\begin{align*} 
\hom_\EuA(a,a)& = CF^*(a^1, a) = \BbK e_0 \oplus \BbK e_1 ;  \\
\hom_\EuA(b,b) &= CF^*(b^1, b) = \BbK f_0 \oplus \BbK f_1,  \end{align*}
for generators $e_0$ and $f_0$ of degree $0$ and $e_1$ and $f_1$ of degree $1$.

We specify the product $\mu^2_\EuA$ by giving instead the structure of the cohomological category $A=H^*(\EuA)$. This graded-linear category has hom-spaces $\Hom(X,Y):=H^* \hom_\EuA(X,Y)$ (which equals $\hom_{\EuA}(X,Y)$ by minimality)  and composition $[x]\circ  [y] = (-1)^{|y|} [\mu^2_\EuA(x,y)]$.  

In $A$, $[e_0]$ (resp. $[f_0]$) is the identity morphism of $\Hom(a,a)$ (resp. $\Hom(b,b)$). The non-zero products  that do not involve $[e_0]$ or $[f_0]$ are 
\begin{equation} [u_1][v_0]=[f_1]; \quad [v_0][u_1]=[e_1] .\end{equation}

We view $A$ as an associative algebra by taking the direct sum of all four morphism spaces; multiplication in this algebra is given by composition when that makes sense, and by zero in other cases. To remember the categorical structure, one regards $A$ as an $(S,S)$-bimodule, where $S$ is the semisimple ring $\mathbb{K}^2=\BbK\{ [e_0],[f_0] \}$. For example, $\Hom_A(a,b)=[f_0] A [e_0]$.

\paragraph{Quiver presentation.} If one views $A$ as a $\BbK$-algebra, one can present it as a quotient of the path-algebra for the following quiver $Q$:
\begin{equation}
\xymatrix{
 a \ar@/^/^{[u_1]}[r] & b\ar@/^/^{[v_0]}[l]
}
\end{equation}
The path-algebra $\Pi Q$ is the unital, graded, associative $\mathbb{K}$-algebra with a basis given by all composable sequences of arrows $a \to b$ or $b \to a$ (written as $[u_1]$ and $[v_0]$ and given degrees 1 and 0 respectively), including the trivial sequences starting and ending at $a$ or at $b$. Multiplication is concatenation when this is results in a composable sequence, and is zero otherwise. The $\BbK$-algebra $\Pi Q$ is actually an $(S,S)$-bimodule:  $[e_0]$ and $[f_0]$ act by composing (on the left or right) with the trivial sequence at $a$ or $b$ respectively.  To obtain $A$ from $\Pi Q$, we divide by the ideal generated by the sequence $a\to b\to a \to b$ and $b\to a\to b\to a$. This ideal is the third power $I^3$ of the `arrow ideal' $I$---the ideal generated by all paths of positive length. There is a canonical isomorphism of graded $\BbK$-algebras $\Pi Q/I^3 \to A$, which sends the arrow $[u_1]$ to the class $[u_1]\in A$, and the arrow $[v_1]$ to $[v_1]\in A$. This isomorphism respects the $(S,S)$-bimodule structure.
\\
\begin{thm}\label{bigradedHH}
Work over an arbitrary base field $\BbK$.

(a) The bigraded Hochschild cohomology 
\[ HH^{r+s}(A,A)^s = \Ext^r_{(A,A)}(A,A[s]) \]
satisfies the relation
\[  HH^{r+s+2}(A,A)^{s-6} \cong HH^{r+s}(A,A)^s ,\quad r>0, \; s\in \Z, \]
which amounts to 8-periodicity in $r$. If $\Char(\BbK)=2$ then one has 
$HH^{r+s+1}(A,A)^{s-3} \cong HH^{r+s}(A,A)^s. $

(b) For $0\leq r\leq 8$, $HH^{r+s}(A,A)^{s}$ is as follows:
\[\begin{array}{c|| c|c|c|c|c|c|c|c|c|c}
 & 0 & 1 & 2 & 3 & 4 & 5 & 6 & 7 & 8 & r \\ \hline\hline
1 & \BbK^2&&&&& &&&\\ \hline
0 & \BbK 	&\BbK &&&&&&&\\ \hline
-1 & & &   \BbK/(2)  &\BbK/(2) &&&&&\\ \hline
-2 &&&&\BbK/(3)&\BbK/(3)&&&&\\ \hline
-3  & &  & & &\BbK/(2) &\BbK/(2)&&&\\ \hline
-4  & & & & & & & \BbK &  \BbK & \\\hline
-5   &&&&&&&&& \\ \hline
-6& & & & & &&  &&\BbK  \\ \hline
s &&&&&&&&&\\ 
\end{array}\]
\end{thm}

The table in (b) is the result of a computation involving rather delicate sign considerations. We are grateful to Paul Seidel, who has written a Python 3.0 program to compute $HH^{r+s}(A,A)^s$ for small $r$, over a field of small positive characteristic, via the normalized bar complex. We have used this program to check our result.

\paragraph{Graded Ext.} The Ext-modules in the statement of the theorem are defined, for graded $(A,A)$-bimodules $B$, as right-derived functors: $\Ext^r_{(A,A)}(A,B)=\mathbf{R}^r\Hom^0_{(A,A)}(A,-)(B)$. Here $\Hom^0_{(A,A)}$ refers to degree-preserving homomorphisms of $(A,A)$-bimodules. Concretely, one takes a projective resolution $P_* \to B \to 0$ by graded $(A,A)$-bimodules (the normalized bar resolution is one possibility).  The resolution is a sequence of bimodules, $P_*=\{\dots \to P_2 \stackrel{p_2}{\to} P_1 \stackrel{p_1}{\to} P_0\}$, and the $p_i$ are degree-preserving bimodule maps.  Applying the Hom-functor results in a cochain complex $\Hom(P_*,B)$, whose maps $p_k^*\colon \Hom(P_{k-1}, B)\to \Hom(P_k, B)$ are given by $p_k^*\theta = \theta\circ p_k$, and one has $\Ext^r(A,B) = H^r \Hom(P_*,B)$.

Shifting a bimodule $B$ changes its structure as a left $A$-module: the multiplication $A\otimes B[s]\otimes A\to B[s]$ is $a\otimes b \otimes a' \mapsto (-1)^{s |a|} aba'$, where $aba'$ is the result of multiplication $A\otimes B\otimes B\to B$.  If $B'$ is another graded bimodule, we write $\Hom^*(B',B)$ for $\bigoplus_s{\Hom^0(B',B[s])}$.

\begin{pf}[Proof of the theorem.]
Sk\"oldberg \cite{Sko} finds a projective resolution of $A$ as a graded $(A,A)$-bimodule. To describe it, let $B_{2j} \subset \Pi Q$ be the subspace of the path-algebra spanned by paths of length $3j$, and let 
$B_{2j+1}\subset \Pi Q$ be the subspace spanned by paths of length $3j+1$.
Then $B_j$ is an $(S,S)$-bimodule, and is a rank 2 free $\BbK$-module. Ordered bases $(\beta_j,\gamma_j)$ for $B_j$ can be given as follows:
\begin{align*}
& B_0 = S = \BbK\{ [e_0],[f_0]\} ; \\
& B_{4k}  = \BbK \{ ( [u_1][v_0])^{3k} , ([v_0] [u_1])^{3k}\} \quad (k>0); \\
& B_{4k+1} = \BbK \{ [v_0]( [u_1][v_0])^{3k} ,  [u_1]([v_0][ u_1])^{3k}\}; \\
& B_{4k+2} = \BbK \{ [v_0] ([u_1][v_0])^{3k+1},[u_1] ([v_0] [u_1])^{3k+1} \}; \\
& B_{4k+3} = \BbK \{ ([u_1][v_0])^{3k+2}, ([v_0] [u_1])^{3k+2} \} .
\end{align*}
So $\beta_0=[e_0]$, etc. The resolution is
\[    
0 \leftarrow A \stackrel{\varepsilon}{\longleftarrow} P_0 
 \stackrel{p_1}{\longleftarrow} P_1
   \stackrel{p_2}{\longleftarrow} P_2
\stackrel{3_2}{\longleftarrow} P_2  \longleftarrow \dots
 \]   
where 
\begin{equation}  P_j =  A\otimes_S B_j \otimes_S A, \end{equation}
$\epsilon(a\otimes b)=ab$, and the maps $p_j\in \Hom^0_{(A,A)}(P_j,P_{j-1})$ are as follows.
Write a typical monomial $a\otimes c_k \otimes\dots \otimes c_1 \otimes b$ as 
$a[c_k | \dots | c_1]b$. Then
\begin{align}
p_{2j}( a[c_{3j} | \dots | c_1]b)  &=  
	a[c_{3j} | \dots | c_3]c_2c_1 b + ac_{3j} [c_{3j-1} | \dots | c_2] c_1 a  \\
	&+ a c_{3j} c_{3j-1} [c_{3j-2} | \dots | c_1] b; \notag\\
p_{2j+1}(a[c_{3j+1} | \dots | c_1]b) &=  
ac_{3j+1}[c_{3j}| \dots | c_1] b -
a [c_{3j} | \dots | c_2] c_1 b .
\end{align}
It's easy to check that $p_k\circ p_{k-1} =0$ and that $\varepsilon\circ p_1=0$. It is shown in \cite{Sko} that this complex is a projective bimodule resolution of $A$.

Observe that $P_j$ is generated, as an $(A,A)$-bimodule, by its subspace $G_j = S \otimes_S B_j \otimes_S S$. One has a basis for $G_{4k}$ given by
$f_0 \otimes \beta_{4k}\otimes f_0$ and $e_0\otimes\gamma_{4k}\otimes e_0$, and hence an isomorphism of graded $\BbK$-vector spaces 
\begin{equation} G_{4k} \cong \BbK^2[-3k].  \end{equation}
Similarly,
\begin{align*}
& G_{4k+1} \cong \BbK[-3k] \oplus \BbK[-3k-1]\\
& G_{4k+2} \cong \BbK[-3k-1] \oplus \BbK[-3k-2]\\
& G_{4k+3} \cong \BbK^2[-3k-2].
\end{align*}
In each case, the isomorphism comes from the ordered basis $(1\otimes \beta_j \otimes 1, 1\otimes \gamma_j \otimes 1)$ for $G_j$. Moreover, $P_j$ is the direct sum  of its sub-bimodules generated by $\beta_j$ and $\gamma_j$.  These summands are not free; one has, for instance, an isomorphism
\[ A\beta_{4k}A \to (Af_0 \otimes f_0 A) [-3k], \quad a \beta_{4k} a'\mapsto af_0\otimes f_0 a'. \]
Restriction to $G_j$ gives an injective map of graded vector spaces
\[ r_j \colon \Hom^*_{(A,A)}(P_j, A)  \to \Hom_\BbK^*(G_j, A),\quad r_j(\lambda)=\lambda|_{G_j} ,  \]
and one has, for example, 
\[ \im(r_{4k})= f_0Af_0 [-3k]  \oplus e_0Ae_0 [-3k].\]
We now consider the maps $p_j^*$. As a sample, consider $p_{4k+1}^*$. If $\theta\in \Hom^s_{(A,A)}(P_{4k},A)$, then
\begin{align*}  
(p_{4k+1}^*\theta)(\beta_{4k+1}) &= (\theta\circ p_{4k+1})(\beta_{4k+1})\\
&=\theta \big( v_0 f_0 [( u_1v_0)^{3k}]f_0 - f_0[ (v_0u_1)^{3k}] v_0 f_0 \big)\\
&=\theta \big( v_0 \beta_{4k} -  \gamma_{4k}v_0 f_0 \big)\\
&=\left (v_0  \theta ( \beta_{4j})  -  \theta(\gamma_{4j})v_0\right) .
\end{align*}

Using this and similar calculations, we can describe the cochain complex $\Hom^*_{(A,A)}(P_*,A)$ in explicit terms. The segment beginning at
\begin{equation} \Hom^*_{(A,A)}(P_{4k},A)
\stackrel{p_{4k+1}^*}{\longrightarrow} \Hom^*_{(A,A)}(P_{4k+1},A)  \end{equation}
reads as follows:
\\

$
\begin{CD}
   f_0Af_0[-3k]   \oplus   e_0Ae_0[-3k]  \\
  @V{ p_{4k+1}^*}V{\left[ \begin{array}{cc}
       l(v) &  -r(v) \\ 
       -r(u)  & \eta l(u)   \\ 
   \end{array}\right] }V\\
   e_0Af_0[-3k]  \oplus  f_0Ae_0[-3k-1] \\ 
   @V{p_{4k+2}^*}V{\left[ \begin{array}{cc}
       \eta l(e_1)+ r(f_1) &  l(v)r(v)    \\ 
       \eta l(u)r(u) &  -\eta l(f_1)+ r(e_1)  \\ 
   \end{array}\right]}V\\
  e_0Af_0[-3k-1]  \oplus  f_0 A e_0[-3k-2] \\ 
  @V{p_{4k+3}^*}V{ \left[ \begin{array}{cc}
      -\eta l(u) & -r(v) \\ 
     -r(u)  &  l(v)  \\ 
  \end{array}\right]}V\\
  f_0Af_0[-3k-2]  \oplus  e_0Ae_0[-3k-2] \\ 
  @V{p_{4k+4}^*}V{ \left[ \begin{array}{cc}
     \eta l(f_1)+ r(f_1) &  \eta l(u) r(v) \\ 
     l(v)r(u)  &  \eta l(e_1)+r(e_1) \\ 
  \end{array}\right]}V\\
   f_0 Af_0[-3k-3]  \oplus e_0 Ae_0 [-3k-3]. \\ 
\end{CD}$
\vspace{3mm}

Here $l(a)(x)=a x$ and $r(b)(x)=x b$ (this is multiplication in $A$, disregarding degree shifts); and $\eta= (-1)^{k}$. These four adjacent maps depend on $k$ only through the alternating sign $\eta$ and the degree shifts. So we see that $H^{r+8} (\Hom(P_*,A[s]) \cong H^r(\Hom(P_*, A[s-6])$ for $r>0$; and that in characteristic 2, $H^{r+4} (\Hom(P_*,A[s]) \cong H^r(\Hom(P_*, A[s-3])$. With this explicit description of the complex, the cohomology groups can be determined by a routine computation which is left to interested readers to reproduce.
\end{pf}

\paragraph{Parametrizing $A_\infty$-structures.}
If we avoid characteristics dividing $6$, the structure of $HH^*(A,A)$ is simple enough that we can use our result to determine the gauge-equivalence classes of minimal $A_\infty$-structures on $A$  (for gauge-equivalence, alias formal diffeomorphism, see \cite{SeiQuartic}.) By homological perturbation theory, these classes correspond naturally to pairs of an $A_\infty$-algebra $\EuScript{B}$ with an isomorphism $H^0\EuScript{B}\to A$, up to quasi-isomorphism.
\\
\begin{thm}\label{classification}
Assume $\Char(\BbK)\notin \{2,3\}$. Consider an arbitrary $A_\infty$-structure $\EuA'$ on $A$, with $\mu^1_{\EuA'}=0$ and $\mu^2_{\EuA'}=\mu^2_\EuA$.  
\begin{enumerate}
\item
The structure $\EuA'$ is equivalent, under the action of the group of gauge transformations, to an $A_\infty$-structure $\EuA''$ in which the structure maps $\mu^1$, $\mu^3$, $\mu^4$ and $\mu^5$ all vanish. 
\item For $\EuA''$ as in (1), the products $\mu^6_{\EuA''}$ and $\mu^8_{\EuA''}$ are cocycles representing classes $m^6\in HH^2(A,A)^{-4}$ and $m^8\in HH^2(A,A)^{-6}$ respectively. These classes are invariants of $\EuA'$. Two $A_\infty$-structures which define the same classes $(m_6,m_8)$ are gauge-equivalent.
\item
Every class in $HH^2(A,A)^{-4} \times HH^2(A,A)^{-6}$ is realized as the invariant $(m^6,m^8)$ of some minimal $A_\infty$-structure on $A$. Hence gauge-equivalence classes of minimal $A_\infty$-structures on $A$ are canonically parametrized  by this product of two 1-dimensional $\BbK$-vector spaces, and the action of $\BbK^*$ rescaling structure maps, defined by $\mu^d_{t\EuA'}= t^{d-2}\mu^d_{\EuA'}$, corresponds to the action $t\cdot(a,b)=(t^4 a , t^6 b)$ on $HH^2(A,A)^{-4} \times HH^2(A,A)^{-6}$.
\end{enumerate}
\end{thm}

\begin{pf} 
In general, if  $\mu^\bullet$ and $\tilde{\mu}^\bullet$ are the structure maps for $A_\infty$-structures on $A$, with $\mu^1=\tilde{\mu}^1=0$ and $\mu^2(x,y)=\mu^2(x,y)=(-1)^{|y|} xy$, and if $\mu^k=\tilde{\mu}^k$ for $k<d$, then the difference $\mu^d- \tilde{\mu}^d \in CC^2(A,A)^{2-d}$ is a cocycle: $\delta (\mu^d- \tilde{\mu}^d)=0$. If it is a coboundary, say $\mu^d-\tilde{\mu}^d=\delta \nu$, then the gauge transformation $G$, with components $g^k\in CC^1(A,A)^{1-k}$ given by $g^1=\id$, $g^{d-1}=\nu$ and $g^k=0$ for $k\notin \{ 1 ,d-1\}$, has the property that $(G_*\tilde{\mu})^k=\mu^k$ for $k\leq d$.

Since $HH^2(A,A)^{2-r}=0$ for $r>8$ by Theorem \ref{bigradedHH}, to prove gauge-equivalence of $\EuA'$ and $\EuA''$ it suffices to make them agree up to $\mu^8$. By taking $\tilde{\mu}^\bullet$ to be the trivial $A_\infty$-structure, and noting that $HH^2(A,A)^{2-r}=0$ for $r\in \{3,4,5\}$, we see that one can find a gauge transformation $g$ so that the $A_\infty$-structure $\mu'^\bullet:=(G_* \mu)^\bullet$ has vanishing structure maps of these orders. One then has $\delta \mu'^6=0$. 
We shall see presently that the class $m^6(\EuA')=[\mu'^6]\in HH^2(A,A)^{-4}$ is an invariant of $\EuA$. 

The $A_\infty$-equation, and the vanishing of $\mu'^r$ for $r\in \{3,4,5\}$, imply that $\mu'^7$ and $\mu'^8$ are again cocycles, representing classes $m^7$ and $m^8$. Now $m^7 \in HH^2(A,A)^{2-7}$, and this space vanishes by Theorem \ref{bigradedHH} (b), so there is a gauge transformation $L$ such that $(L_*\mu')^7=0$ while $(L_*\mu')^d=\mu'^d$ for $d<7$. Hence two $A_\infty$-structures which have the same class $m^6$ are gauge-equivalent up to order 7. If in addition their $m^8$-classes agree then a suitable gauge transformation makes their structure maps agree up to order 8, whereupon the structures are gauge-equivalent. 

There exists a minimal $A_\infty$-structure with prescribed values for $(m^6,m^8)$. Begin by declaring $\mu^k=0$ for $k=3,4,5$. If $\mu^k$ has been chosen for $k<d$, and the $A_\infty$ relations hold for all input-sequences of length $< d$, then the condition on $\mu^d$ for the $A_\infty$ relations to hold on sequences of length $d$ is that
\begin{equation}  \label{Maurer Cartan}
\delta \mu^d = \frac{1}{2}\sum_{j=3}^{d-1} {[\mu^{j}, \mu^{d-j+2}]}. 
\end{equation}
Here $[\cdot , \cdot]$ is the Gerstenhaber Lie bracket \cite[(3b)]{SeiQuartic}.
By applying the same identity for $j$ instead of $d$, and the Jacobi identity, one sees that the right-hand side defines a cocycle in $CC^3(A,A)^{3-(d+1)}$, and to define $\mu^d$ it suffices to check that this cocycle is a coboundary. But $HH^3(A,A)^{3-(d+1)}=0$ except for $d\in \{ 6,8\}$. For those values of $d$, we take $\mu^d$ to be a Hochschild cocycle representing the given cohomology class (so $\delta \mu^6=\delta\mu^8=0$), which is legitimate since the Gerstenhaber brackets on the right-hand side of (\ref{Maurer Cartan}) are all zero in these cases. Inductively, one can define $\mu^d$ for all orders $d$.

The remaining question is the independence of $m^6$ and $m^8$ from the gauge transformations used to define them.   Here we appeal to indirect arguments. The well-definedness of $m^6$ can be seen as a consequence of its interpretation as a differential in the length spectral sequence (see point $\clubsuit$ below), which is itself invariant under gauge transformations. We can see that both $m^6$ and $m^8$ are invariant by exhibiting a 2-parameter family of non-gauge-equivalent $A_\infty$-structures which realize all possible values of $(m^6,m^8)$. Ext-algebras on a family of Weierstrass cubics will achieve this (see $\spadesuit$ below).
 \end{pf}

\paragraph{The length spectral sequence.}  $A_\infty$-categories $\EuScript{B}$ also have Hochschild cohomology $HH^*(\EuScript{B},\EuScript{B})$. This graded $\BbK$-vector space is not bigraded, but the underlying bar complex $CC^*(\EuScript{B},\EuScript{B})$ has a natural filtration, the \emph{length filtration}. The resulting spectral sequence $\{( E_k^{rs},d^k_{rs}) \}$ has $E_2^{rs}= HH^{r+s}(B,B)^s$, where $B=H^*(\EuScript{B})$ \cite[(1f)]{SeiBook}.  In general, convergence of the spectral sequence is not guaranteed, but from Theorem \ref{bigradedHH} we obtain the following two corollaries:
\\
\begin{cor}
For any $A_\infty$-structure $\EuA'$ on $A$, with $\mu^1_{\EuA'}=0$ and $\mu^2_{\EuA'}(x,y)=(-1)^{|y|}xy$, the length spectral sequence $\{ (E^{rs}_k,d_k^{rs}) \}$ for the length filtration on $CC^*(\EuA',\EuA')$ degenerates by $E_{10}$.
\end{cor}
\begin{pf}
One has $E_2^{rs}=HH^{r+s}(A,A)^s$. The longest possible differentials in the  spectral sequence are of form $d_9^{8k+7,8k-4}$. 
\end{pf}

\begin{cor}
Assume $\Char(\BbK)\notin \{ 2,3 \}$. Let $\EuA'$ be a minimal $A_\infty$-structure on $A$.  Then either $\EuA'$ is formal, or $d_5\neq 0$, or $d_7 \neq 0$ in the length spectral sequence for $HH^*(\EuA',\EuA')$.  The only possibly non-trivial differentials in the length spectral sequence  are $d_5$, $d_7$ and $d_9$.  One has 
\begin{align*} 
& \dim HH^0(\EuA',\EuA') =1; \\
&\dim HH^1(\EuA',\EuA')\leq 3; \\
&\dim HH^d(\EuA',\EuA') \leq 2 \; (d>1). 
\end{align*}
If $\EuA'$ is not formal then $\dim HH^1(\EuA',\EuA')=2$ and $\dim HH^2(\EuA',\EuA')=1$.
\end{cor}
\begin{pf}
In the spectral sequence, one has $d_j=0$ for $2\leq j\leq 4$, while $d_5$ is the Gerstenhaber bracket $[m^6,\cdot]$ \cite[(3.7)]{SeiQuartic}. If $m^6=0$ then $d_5=d_6=0$ and $d_7=[m^8,\cdot]$. The Euler derivation $e \colon A\to A$, given by $e(x)=\deg(x) x$ for homogeneous $x$, represents a non-trivial class $[e]\in HH^1(A,A)^0\cong \BbK$ and one has $d_{k-1}[e] =(k-2) m^k$ when $d_j=0$ for $2\leq j\leq k-2 $ \cite[(3.8)]{SeiQuartic}.  ($\clubsuit$ This explains why $m^6$ is well-defined.) Armed with this information, and observing that the  identity natural transformation $\id_{\EuA'}\to \id_{\EuA'}$ represents a non-trivial class in $HH^0(\EuA',\EuA')$, the assertions are straightforward to check.
\end{pf}
\paragraph{Non-formality.} One can compare this corollary to results and conjectures in symplectic geometry. A conjecture of Seidel's \cite{SeiDeform} states that, for certain Liouville domains $M$, the `open-closed string map' 
\begin{equation} \kappa_M \colon  SH^*(M)\to HH^*(\EuF(M),\EuF(M))  \end{equation}
from symplectic cohomology is an isomorphism.  According to \cite{SeiBias}, one has
\begin{equation} SH^0(T_0)=\BbK,\; SH^1(T_0) = \BbK^2,\; SH^d(T_0)=\BbK \; (d>1).
\end{equation}
There is an isomorphism $HH^*(\EuF(T_0),\EuF(T_0))\to HH^*(\EuA,\EuA)$, so the conjecture predicts that $\EuA$ is non-formal. This is true:
\\
\begin{thm}\label{m6}
When 6 is invertible, the $A_\infty$-structure $\EuA$ arising from $\wh{\EuF}(T_0)$ is not formal. Indeed, $m^6(\EuA)\neq 0$.
\end{thm}
Our proof of this is computational: we take a dg model for $\EuA$ due to Abouzaid \cite{Abo}, apply the homological perturbation lemma to transfer it to an $A_\infty$-structure on $A$, then explicitly find a gauge transformation that kills $\mu^3$ and $\mu^4$. The resulting $\mu^6$ can then be seen to represent a non-zero class $m^6(\EuA)$. The details of our computer-assisted calculation are given in the Appendix \ref{appendix}.

\begin{rmk}
One can also consider the closed torus $T$. There is again an open-closed string map, which is an isomorphism \cite{AS} of graded $\Lambda_\C$-algebras
\begin{equation} \kappa_T \colon  H^*(T; \Lambda_{\C}) \to HH^*(\EuA(T),\EuA(T)). \end{equation} 
The base field (on both sides) is the complex Novikov field $\Lambda_{\C}$, consisting of formal series $\sum_{r\in \R}{a(r) t^r}$,  where $a\colon \R\to \C$ is a function such that $\mathrm{supp}\, (a)\cap (-\infty,c] $ is finite for any $c\in \R$.  On the right-hand side, $\EuA(T)$ is the full subcategory of the Fukaya category of $T$ with objects $a$ and $b$. 
\end{rmk}

\paragraph{Weierstrass cubics.} We draw the reader's attention to an algebro-geometric interpretation of $HH^2(A,A)$. We have seen that when $1/6\in \BbK$, gauge-equivalence classes of $A_\infty$-structures on $A$ are parametrized by the 2-dimensional vector space $HH^2(A,A)$, on which the $\BbK^*$-action  rescaling the structure maps has weights $(4,6)$. This is reminiscent of the family of cubic curves $\EuScript{E}\subset \mathbb{A}^2 \times \mathbb{P}^2$ defined by the the Weierstrass equation $y^2 z =4x^3-pxz^2-qz^3$ considered as a family over $\mathbb{A}^2$ (see e.g.  \cite[ch. 2]{KM}). Here $(p,q)$ are coordinates on $\mathbb{A}^2$, $[x:y:z]$ homogeneous coordinates on $ \mathbb{P}^2$.  Each fiber in this family carries a non-zero differential $dx/(2y)$, and therefore by Serre duality a basis for $H^1(\mathcal{O})$. Over the locus where the discriminant $\Delta:=p^3 - 27q^2$ is non-zero, the fiber $\EuScript{E}_{(p,q)}$ is an elliptic curve.  When $\Delta=0$ and $(p,q)\neq (0,0)$, $\EuScript{E}_{(p,q)}$ has a node; and when $(p,q)=(0,0)$, it has a cusp. The action of $\BbK^*$ on the base $\BbK^2=\mathbb{A}^2(\mathbb{K})$, $t\cdot(p,q)=(t^4 p,t^6 q)$, is covered by an action on the family, $t\cdot (x,y,z)=(x,ty,t^{-2}z)$.
 
There are two distinguished objects in the dg category of coherent sheaves on $\EuScript{E}_{(p,q)}$: the structure sheaf $\mathcal{O}$, and the skyscraper sheaf $\mathcal{S}$ at the point $P=[0:1:0]$ at infinity, which is always a regular point. The full subcategory $\EuScript{B}_{(p,q)}$ generated by these two objects comes with an isomorphism $A\to H\EuScript{B}_{(p,q)}$. One has $\Ext^*(\mathcal{O},\mathcal{O})\cong \Lambda^*H^1(\mathcal{O})$, but we can use the basis for $H^1(\mathcal{O})$ to identify this exterior algebra with $f_0Af_0$. This explains why $A\to H\EuScript{B}_{(p,q)}$ is canonical. 

\begin{prop}\label{Weierstrass}
The map $M\colon \BbK^2 \to HH^2(A,A)$, $(p,q)\mapsto m^6(\EuScript{B}_{(p,q)})+m^8(\EuScript{B}_{(p,q)})$, is a linear bijection.
\end{prop} 
\begin{pf}
The crucial observation is that   one can rebuild the homogeneous coordinate ring $\bigoplus_{n\geq 0}{H^0(\EuScript{E}_{(p,q)}, \mathcal{O}(nP))}$---and hence the elliptic curve $\EuScript{E}_{(p,q)}$---from the dg category $\EuScript{B}_{(p,q)}$, or equally from its $A_\infty$ minimal model. One inductively expresses $\mathcal{O}(nP)$ as a twisted complex in $\mathcal{O}$ and $\mathcal{S}$ via exact triangles $\mathcal{O}((n-1)P) \to\mathcal{O}(nP) \to \mathcal{O} \to \mathcal{O}((n-1)P)[1]$ and $\mathcal{O}\to \mathcal{O}(1) \to \mathcal{S}\to \mathcal{O}[1]$. The ring structure is then realized via composition in $\tw \EuScript{B}_{(p,q)}$. One also recovers the basis of $H^1(\EuScript{E}_{(p,q)},\mathcal{O})$ Serre-dual to the differential $dx/(2y)$. These data are sufficient to determine $(p,q)$. This proves that $(m^6,m^8)$ are gauge-invariants of an $A_\infty$-structure ($\spadesuit$), and also shows that $M$ is bijective.

We next claim that $M$ is a polynomial map (we thank Paul Seidel for technical assistance with this point).  Being polynomial and $\BbK^*$-equivariant, $M$ is necessarily linear.

The dga $(\EuScript{B},d_{\EuScript{B}}^*)$ can be realized, over the base $R:=\BbK[p,q]$, as $\mathrm{End}\,(C^*)$, where the cochain complex $C^*$ of $R$-modules $C^*$ is the \v{C}ech complex for an open cover by two affines, with coefficients in the sheaf $\mathcal{E}nd(\mathcal{O}\oplus \mathcal{S})$. One checks that $C^*$ is a complex of free $R$-modules.  The cohomology $H^*\EuScript{B}=A\otimes R$ is a free $R$-module, so one can write $\ker d_{\EuScript{B}}^* = \EuScript{H} \oplus \im d_{\EuScript{B}}$ for a submodule $\EuScript{H}$ projecting  isomorphically  to $H^*\EuScript{B}$. Since $\mathrm{End}(C^*)$ is a projective $R$-module, we can split $d_{\EuScript{B}}^*\colon \mathrm{End}(C^*)\to \im d_{\EuScript{B}}^*$ and so write $\mathrm{End}(C^*)= \ker d_{\EuScript{B}} \oplus \EuScript{C} =  \im d_{\EuScript{B}}^* \oplus \EuScript{H} \oplus \EuScript{C}$. The projection map $\EuScript{B}\to \EuScript{H}$ and the inclusion $\EuScript{H}\to \EuScript{B}$ are then homotopy-inverses; there is a contraction of $\EuScript{B}$ to $A\otimes R$. Applying the homological perturbation lemma to this contraction, one obtains an $A_\infty$-structure on $A\otimes R$, over $R$. Theorem \ref{bigradedHH} remains valid over $R$ (a flat $\BbK$-module). Points 1-2 of Theorem \ref{classification} are also valid---the proof involves no division---and so one obtains a class $(m^6,m^8)\in HH^2_R(A\otimes R,A\otimes R)$. The latter module is isomorphic to $HH^2_{\BbK}(A)[p,q]$ by a canonical  map sending $(m^6,m^8)$ to $M$; hence $M$ is polynomial. 
\end{pf}
A natural follow-up will be to compare this result with Polishchuk's calculation \cite{Pol}.   We will explain elsewhere that $\EuA\subset \wh{\EuF}(T_0)$ has parameters $(p,q)=M^{-1}(m^6,m^8)$ which lie on the discriminant curve $\{\Delta=0\}$. Since $p\neq 0$, it follows that $\EuF(T_0)$ is quasi-isomorphic to the category $\EuScript{B}_{(p,q)}$ for a nodal cubic---an instance of homological mirror symmetry.

{\bf Acknowledgments:} YL thanks MIT for generous hospitality during the writing of this paper. TP was partly supported by  by NSF grant 1049313.  We both thank Mohammed Abouzaid, David Ben-Zvi, Pavel Etingof and Paul Seidel for sharing their ideas and expertise, and we also thank Seidel for writing a computer program that provided a check on our calculations.

\newpage
\appendix
\section{Appendix: The Fukaya $A_\infty$-algebra of two plumbed circles is not formal}
\label{appendix}

Let $\EuA$ be the $A_\infty$-algebra given by the direct sum of
the four hom-spaces in the full subcategory of $\EuF(T_0)$---the exact Fukaya
category of a once-punctured 2-torus---whose objects are the two oriented
Lagrangian branes. The underlying curves meet transversely at a single point,
and their fundamental classes span $H_1(T_0)$. They are equipped with gradings
such that all hom-spaces are supported in degrees 0 and 1, and with non-trivial
spin-structures.

We showed that, over a field $\BbK$ in which 6 is invertible, $\EuA$
determines a primary deformation class $m^6(\EuA)\in HH^2(A,A)^{2-6}$, where
$A$ is the $\BbK$-algebra $H^*\EuA$. This class is a quasi-isomorphism
invariant of $\EuA$; if it is non-zero, $\EuA$ is necessarily non-formal, i.e.,
not quasi-isomorphic to $A$ with the `trivial' $A_\infty$-structure (the one
with all structure-maps vanishing except the 2-fold product). We have verified,
with the assistance of some simple computer programs written in Python 3.0,
that \[m^6(\EuA)\neq 0.\] In this appendix we document this calculation. 

\paragraph{DG models.} Consider a symplectic manifold with two embedded
Lagrangian submanifolds $Q_1$ and $Q_2$ intersecting transversely. Let $M$ be a
regular neighborhood of $Q_1\cup Q_2$.  Abouzaid's article \cite{Abo} describes
(under certain assumptions) a dg algebra quasi-isomorphic to the full
subcategory of the exact Fukaya category of $M$ whose objects are $Q_1$ and
$Q_2$ equipped with gradings and relative spin structures.  

We apply this general result to obtain the following dg category 
$\EuScript{D}$, with objects $a$ and $b$, quasi-isomorphic to $\EuA$:

\begin{center}
  \begin{tabular}{@{} c|cc}
   hom-space & basis for $\hom^0$ & basis for $\hom^1$\\ 
\hline   
 $\hom_{\EuScript{D}}(a,a)$ & $x_0$, $x_1$, $x_2$ & $x_{01}$, $x_{12}$, $x_{02}$ \\ 
 $\hom_{\EuScript{D}}(b,b)$ & $y_0$, $y_1$, $y_2$ & $y_{01}$, $y_{12}$,$y_{02}$ \\ 
 $\hom_{\EuScript{D}}(b,a)$ & $v_0$, $v_1$ &  $v_{01}$\\
 $\hom_{\EuScript{D}}(a,b)$ && $u_{01}$\\
  \end{tabular}
\end{center} 
The differential $d \colon \EuD^0 \to \EuD^1$ is as follows:
\begin{align*}
&d x_0 = -x_{01} -x_{02}, 
&& d x_1 = x_{01} - x_{12},
&&& d x_2 = x_{12} + x_{02}, \\
&d y_0 = -y_{01} -y_{02},
&& d y_1 = y_{01} - y_{12},
&&&d y_2 = y_{12} + y_{02} ,\\
&d v_0 = -v_{01}, 
&& d v_1 = v_{01}. 
\end{align*}
The product is as follows. We use the $A_\infty$ (not associative) sign convention; to obtain a genuine dg category, put $xy=(-1)^{|y|}\mu^2_{\EuD}(x,y)$.

\[ 
\begin{array}{l | l}
\hom(a,a) \otimes \hom(a,a) \to \hom(a,a) & \hom(b,b) \otimes \hom(b,b) \to \hom(b,b)\\
\hline
\mu^2_{\EuScript{D}}(x_0, x_0) = x_0 & \mu^2_{\EuScript{D}}(y_0, y_0) = y_0  \\
\mu^2_{\EuScript{D}}(x_1, x_1) = x_1 & \mu^2_{\EuScript{D}}(y_1, y_1) = y_1 \\
\mu^2_{\EuScript{D}}(x_2 ,x_2) = x_2 & \mu^2_{\EuScript{D}} (y_2 ,y_2) = y_2 \\ 
\mu^2_{\EuScript{D}}(x_0 , x_{01}) = -x_{01} & \mu^2_{\EuScript{D}}(y_0 ,y_{01}) = -y_{01} \\
\mu^2_{\EuScript{D}}(x_{01}, x_1) = x_{01} &   \mu^2_{\EuScript{D}}(y_{01}, y_1) = y_{01} \\ 
\mu^2_{\EuScript{D}}(x_1 ,x_{12}) = -x_{12} & \mu^2_{\EuScript{D}} (y_1 ,y_{12}) = -y_{12} \\
\mu^2_{\EuScript{D}}(x_{12}, x_2) = x_{12} &   \mu^2_{\EuScript{D}}(y_{12}, y_2) = y_{12} \\
\mu^2_{\EuScript{D}}(x_0 ,x_{02}) = -x_{02}&  \mu^2_{\EuScript{D}} (y_0 ,y_{02}) = -y_{02} \ \\
\mu^2_{\EuScript{D}}(x_{02}, x_2) = x_{02} &  \mu^2_{\EuScript{D}} (y_{02}, y_2) = y_{02}
\end{array}\]

\[ \begin{array}{l | l}
\hom(b,a) \otimes \hom(b,b) \to \hom(b,a) & \hom(a,a) \otimes \hom(b,a) \to \hom(b,a)\\
\hline
\mu^2_{\EuScript{D}}(v_0 ,x_0 )= v_0 & \mu^2_{\EuScript{D}}(y_0 ,v_0 )= v_0, \\
\mu^2_{\EuScript{D}}(v_1, x_1) = v_1 & \mu^2_{\EuScript{D}}(y_1, v_1) = v_1, \\
\mu^2_{\EuScript{D}}(v_{01}, x_1)= v_{01} & \mu^2_{\EuScript{D}} (y_{01}, v_1) = v_{01}\\ 
\mu^2_{\EuScript{D}}(v_0, x_{01})=  -v_{01} &\mu^2_{\EuScript{D}} (y_0 ,v_{01})=  -v_{01}
\end{array}\]

\[ \begin{array}{l | l}
\hom(b,b) \otimes \hom(a,b) \to \hom(a,b) & \hom(a,b) \otimes \hom(a,a) \to \hom(a,a)\\
\hline
\mu^2_{\EuScript{D}}(x_0 ,u_{01})= -u_{01} & \mu^2_{\EuScript{D}}(u_{01}, y_1)=u_{01} \\
\end{array}\]

\[ \begin{array}{l | l}
\hom(b,a) \otimes \hom(a,b) \to \hom(a,a) & \hom(a,b) \otimes \hom(b,a) \to \hom(b,b)\\
\hline
\mu^2_{\EuScript{D}}(v_0, u_{01}) = -y_{01}  & \mu^2_{\EuScript{D}}(u_{01} , v_1 )= x_{01}
\end{array}\]

Define
\[ e_0= x_0+x_1+x_2 \in \hom^0(a,a);\quad f_0 = y_0 +y_1 + y_2 \in \hom^0(b,b); \]
\[  e_1= x_{01}\in \hom^1(a,a); \quad f_1 = y_{01}\in \hom^1(b,b). \]
Let $\EuScript{C}$ be the subcategory of $\EuScript{D}$ with the same objects, and hom-spaces spanned by $e_0$, $f_0$, $e_1$, $f_1$, $v_0$, $v_1$ ,$v_{01}$, $u_{01}$. One checks that this is indeed a subcategory, and that it is preserved by $d$. The inclusion $\EuScript{C}\to \EuScript{D}$ is a quasi-isomorphism. Since it is considerably smaller, we shall henceforth use $\EuScript{C}$ instead of $\EuScript{D}$. Its structure is as follows.

\begin{center}
  \begin{tabular}{@{} c|cc}
   hom-space & basis for $\hom^0$ & basis for $\hom^1$\\ 
\hline   
 $\hom_{\EuScript{C}}(a,a)$ & $e_0$ & $e_1$ \\ 
 $\hom_{\EuScript{C}}(b,b)$ & $f_0$, & $f_1$ \\ 
 $\hom_{\EuScript{C}}(b,a)$ & $v_0$, $v_1$ &  $v_{01}$\\
 $\hom_{\EuScript{C}}(a,b)$ && $u_{01}$\\
\end{tabular}
\end{center} 
The differential $\mu^1_{\EuC}\colon \hom^0_{\EuScript{C}} \to \hom^1_{\EuScript{C}}$ is zero on basis-vectors except that
\[ \mu^1_{\EuC} v_0 = - v_{01}, \quad \mu^1_{\EuC} v_1 = v_{01} . \]
The non-zero products in $\EuScript{C}$ are as follows:
\[ 
\begin{array}{l | l}
\hom(a,a) \otimes \hom(a,a) \to \hom(a,a) & \hom(b,b) \otimes \hom(b,b) \to \hom(b,b)\\
\hline
\mu^2_{\EuScript{C}}(e_0, e_0) = e_0 & \mu^2_{\EuScript{C}}(f_0, f_0) = f_0  \\
\mu^2_{\EuScript{C}}(e_0, e_1) = -e_1 & \mu^2_{\EuScript{C}}(f_0, f_1) = -f_1 \\
\mu^2_{\EuScript{C}}(e_1 ,e_0) = e_1 & \mu^2_{\EuScript{C}}(f_1 ,f_0) = f_1  
\end{array}\]
\[
\begin{array}{l | l}
\hom(b,a) \otimes \hom(b,b) \to \hom(b,a) 	& \hom(a,a) \otimes \hom(b,a) \to \hom(b,a)\\
\hline
\mu^2_{\EuScript{C}}(v_0, f_0) = v_0 & \mu^2_{\EuScript{C}}(e_0, v_0) = v_0  \\
\mu^2_{\EuScript{C}}(v_1, f_0) = v_1 & \mu^2_{\EuScript{C}}(e_0, v_1) = v_1 \\
\mu^2_{\EuScript{C}}(v_{01}, f_0) = v_{01} & \mu^2_{\EuScript{C}}(e_1 ,v_1) = v_{01}  \\
\mu^2_{\EuScript{C}}(v_{0}, f_1) = -v_{01} & \mu^2_{\EuScript{C}}(e_0 ,v_{01}) = -v_{01}  
\end{array}\]

\[
\begin{array}{l | l}
\hom(b,b) \otimes \hom(a,b) \to \hom(a,b) 	& \hom(a,b) \otimes \hom(a,a) \to \hom(a,b)\\
\hline
\mu^2_{\EuScript{C}}(f_0, u_{01}) = -u_{01} & \mu^2_{\EuScript{C}}(u_{01}, e_0) = u_{01}  \\
\end{array}\]

\[
\begin{array}{l | l}
\hom(b,a) \otimes \hom(a,b) \to \hom(a,a) 	& \hom(a,b) \otimes \hom(b,a) \to \hom(b,b)\\
\hline
\mu^2_{\EuScript{C}}(v_0, u_{01}) = -e_1 & \mu^2_{\EuScript{C}}(u_{01}, v_1) = f_1  \\
\end{array}\]

\paragraph{Applying the homological perturbation lemma.}
Observe next that there is a direct sum splitting $\EuScript{C}=\EuScript{G}\oplus \EuScript{H}$, where ${\mu^1_\EuC}|_{\EuScript{H}} =0$,
$\mu^1_\EuC(\EuScript{G})\subset \EuScript{G}$, and $H^*\EuScript{G}=0$. We take $\EuScript{H}$ to be the direct sum of
$\hom_{\EuScript{C}}(a,a)$, $\hom_{\EuScript{C}}(b,b)$, $\hom_{\EuScript{C}}(b,a)$, and the subspace $\BbK\{v_0+v_1\}\subset \hom^0_{\EuScript{C}}(b,a)$, while $\EuScript{G}=\BbK \{v_1, v_{01} \} \subset \hom_{\EuScript{C}}(b,a)$. 

The expression $\EuScript{C}=\EuScript{G}\oplus \EuScript{H}$ defines a projection map $p\colon \EuScript{C}\to \EuScript{H}$, split by the inclusion map $i\colon \EuScript{H}\to \EuScript{C}$ is the inclusion; and a nullhomotopy $T\colon \EuScript{C}\to \EuScript{C}$ of $i\circ p -\id_{\EuScript{C}} $. One has $T(v_{01})=-v_1$, and $Ta=0$ for the remaining basis-vectors $a$. 

We have $A\cong \EuScript{H}$; the isomorphism we shall use is $e_i\leftrightarrow e_i$, $f_i\leftrightarrow f_i$, $u\leftrightarrow u_{01}$, $v\leftrightarrow v_0+v_1$.

The homological perturbation lemma (we use the conventions of \cite[Prop. 1.12]{SeiBook}) determines a minimal $A_\infty$-structure $\EuScript{B}$ on $A$ (minimal means that $\mu^1_\EuB=0$) and an $A_\infty$ quasi-isomorphism $\EuScript{I}\colon \EuScript{B} \to \EuScript{C}$ extending $i=\EuScript{I}^1$. This structure is defined by an explicit recursion: for $d\geq 2$,
\begin{align*}
\EuI^d(a_d,\dots,a_1) &= \sum_{0< m< d}{T\circ \mu^2_{\EuC}\left(\EuI^{d-m}(a_d,\dots, a_{m+1}), \EuI^m(a_m,\dots, a_1)\right)}, \\
\mu^d_\EuB(a_d,\dots,a_m) &= \sum_{0< m< d}{p \circ \mu^2_{\EuC}
\left( \EuI^{d-m}(a_d,\dots, a_{m+1}), \EuI^m(a_m,\dots, a_1)\right)}.
\end{align*} 

\begin{lem}
The only non-zero products $\mu^d_{\EuScript{B}}$ with $d>2$ are 
\[ \mu_{\EuScript{B}}^d ( u, \stackrel{d-3}{\overbrace{e_1,\ldots,e_1}}, v,f_1) = (-1)^{d+1} f_1;\quad 
\mu^d_{\EuScript{B}}( u, e_1,\ldots,e_1, v) = (-1)^d f_1. \]
\end{lem}
\begin{pf}
Since the image of $T$ is $\BbK v_1$, one has $\EuI^d(a_d,\dots,a_1)\in \BbK v_1$. But $\mu^2_{\EuC}(v_1,v_1)=0$, so
\begin{align*} 
\EuI^d(a_d,\dots,a_1) &=
 T \circ \mu^2_{\EuC}
\left(i(a_d), \EuI^{d-1}(a_{d-1},\dots, a_1)\right)+  
T \circ \mu^2_{\EuC}\left (\EuI^{d-1}(a_{d},\dots, a_2),i(a_1) \right),\\
  \mu^d_\EuB(a_d,\dots,a_1) &= p \circ \mu^2_{\EuC}
\left(i(a_d), \EuI^{d-1}(a_{d-1},\dots, a_1)\right)+  
p \circ \mu^2_{\EuC}\left (\EuI^{d-1}(a_{d},\dots, a_2),i(a_1) \right).  
\end{align*}
for $d>2$. Inspecting the multiplication table of $\EuC$, we notice that 
\[ T \circ \mu^2_\EuC(v_1,\cdot)=0,\quad  p\circ \mu^2_\EuC(v_1,\cdot)=0, \]
so in fact
\begin{align*} 
\EuI^d(a_d,\dots,a_1) &=
 T \circ \mu^2_{\EuC}
\left(i(a_d), \EuI^{d-1}(a_{d-1},\dots, a_1)\right),\\
  \mu^d_\EuB(a_d,\dots,a_1) &= p \circ \mu^2_{\EuC}
\left(i(a_d), \EuI^{d-1}(a_{d-1},\dots, a_1)\right),  
\end{align*}
even for $d=2$. By examining $T \circ \mu^2_\EuC(\cdot,v_1)$ we find that the non-zero terms in $\EuI^2$ are
\[  \EuI^2(v,f_1) =  v_1, \quad  \EuI^2(e_1,v) = - v_1,  \]
and that  if  $\EuI^d(a_d,\dots,a_1)\neq 0$ for some sequence of basis vectors $a_j$ with $d>2$ then we must have $a_d=e_1$. Hence, by induction, 
\[  \EuI^d(e_1,\dots,e_1 ,v,f_1) =  (-1)^d v_1, \quad  \EuI^d(e_1,\dots, e_1,v) = (-1)^{d+1} v_1 \quad (d>2) \]
all other outputs for $\EuI^d$ being zero. 

Now assume $d>2$. The map $p\circ \mu^2_\EuC(\cdot, v_1)$ is non-zero only on $u_{01}$, and so if $\mu^d_{\EuB}(a_d,\dots,a_1)\neq 0$ then $a_d=u$. Hence the only non-zero  $\mu^d_{\EuB}$ products are 
\[ \mu^d_\EuB(u,e_1,\dots ,e_1,v,f_1),\quad \mu^d_\EuB (u,e_1,\dots, e_1,v),\] 
and by a straightforward induction, these are as claimed.
\end{pf}

\paragraph{Gauging away $\mu^3$ and $\mu^4$.}
A \emph{gauge transformation} $G$ for $\EuScript{B}$ (also known as a formal
diffeomorphism) is a sequence of degree-preserving $\BbK$-linear maps $g^i\colon
\EuScript{B}^{\otimes d}\to \EuScript{B}[1-d]$, starting with
$g^1=\id_{\EuScript{B}}$. Gauge transformations form a group which acts on the left on
minimal $A_\infty$-structures, defined by the following explicit recursive formulae:

\begin{align*}  
(G_*\mu)^d(a_d,\dots,a_1)
= \sum_{r} \sum_{s_1,\ldots,s_r} ({G_*\mu)^{r}(g^{s_r}(a_d, \ldots, a_{d-s_r+1}),
\ldots, g^{s_1}(a_{s_1}, \ldots, a_1})) \\
+ \sum_{i+j\leq d}{ (-1)^{|a_1|+\ldots+|a_i|+i}g^{d-j+1}(a_d, \ldots, a_{i+j+1},
\mu^{j}(a_{i+j},\ldots,a_{i+1}), a_i\ldots, a_1}).  
 \end{align*}
The sum in the first term is over all $r \geq 2$ and partitions $s_1+\ldots+s_r=d$.
 We write down a gauge transformation $G$ such that $\mu^3_{G_*\EuScript{B}}=0$. Its terms $g^i$ with $i>1$ are zero except for $g^2=:g$, which is as follows:
\begin{align*}
&g(e_1,e_1) = (1/2)e_1, && g(f_1,f_1) = - (1/2) f_1 ; \\
& g(e_1,v) =  - (1/2) v, && g(v,f_1) = (1/2) v ; \\
& g(u,e_1) = - (1/2) u ; && g(f_1,u) = - (1/2) u.
 \end{align*}
Our first computer-aided check is that $G_*\EuScript{B}$ has $\mu^3=0$. The fourth-order product $\mu^4_{G_*\EuScript{B}}$ then has the following non-zero terms:
\begin{align*}
& \mu^4 (e_1,v,f_1,u) = (1/4) e_1, && \mu^4 (e_1,v,u,e_1) = (1/4) e_1, \\ 
& \mu^4 (v,f_1,f_1,u) = -(1/4) e_1, && \mu^4 (v,f_1,u,e_1) = -(1/4) e_1, \\
& \mu^4 (f_1,u,e_1,v) = (1/4) f_1 , && \mu^4 (f_1,u,v,f_1)  = -(1/4) f_1, \\ 
& \mu^4 (u,e_1,v,f_1) = -(1/4) f_1, && \mu^4 (u,v,f_1,f_1) = - (1/2) f_1 ,\\
& \mu^4 (u,e_1,e_1,v) = (3/4) f_1 , && \mu^4 (v,u,e_1,v) = -(1/2) v ,\\ 
& \mu^4 (v,u,v,f_1) = (1/2) v , && \mu^4 (u,e_1,v,u) = (1/2) u ,\\ 
&\mu^4 (u,v,f_1,u) = -(1/2) u.
\end{align*}

To kill $\mu^4$, we now apply another gauge transformation $H$, whose terms $h^i$ with $i>1$ are zero except for $h^3=:h$. 
The trilinear function $h$ is as follows: 

\begin{align*}
&h(v,f_1,u) = - (1/12) e_0 ,
&&h(v,u,e_1) = -(1/12) e_0\\
&h(e_1,e_1,e_1) = (1/3) e_1, 
&&h(f_1,u,v) = -(1/12) f_0 \\
&h(u,e_1,v) = -(1/12) f_0 ,
&&h(f_1,f_1,f_1) = (1/3) f_1 \\
&h(e_1,v,f_1) = - (1/3) v , 
&&h(v,f_1,f_1) = - (1/6) v \\
&h(e_1,e_1,v) = (1/3) v ,
&&h(f_1,f_1,u) = (5/12) u \\
&h(f_1,u,e_1) = (1/3) u ,
&&h(u,e_1,e_1) = (5/12) u 
\end{align*}
The structure $\EuB_{final}:=H_*G_*\EuB$ has vanishing $\mu^3$ and $\mu^4$. (Checking that $\mu^4$ vanishes involves enough equations that computer aid is very useful here.) It is possible to kill $\mu^5$ as well, but we do not need to do so, because the
sixth-order product $\mu^6$ already defines a Hochschild cocycle $\mu^6\in CC^2(A,A)^{2-6}$. 

\paragraph{Non-triviality of the sixth-order product.} We have carried out a computer-aided calculation of $\mu^6$ in $\EuB_{final}$, and verified  that the computer indeed outputs a cocycle. Rather than tabulating the results, we shall instead pick out a small subset of the data and use it to verify that the class $m^6=[\mu^6]\in HH^2(A,A)^{2-6}$ is non-zero. Though computer-guided, this part of the argument is human-checkable.

Suppose that $\mu^6=144 \delta (\nu)$ for some $\nu \in CC^1(A,A)^{-4}$ (the reason for the normalization will become apparent), where $\delta$ is the Hochschild coboundary operator. Then, in $\EuB_{final}$,
\begin{align*} 
-\mu^6(a_6,a_5,a_4,a_3,a_2,a_1)  
=  &- \mu^2(\nu(a_6,a_5,a_4,a_3,a_2),a_1) - \mu^2(a_6,\nu(a_5,a_4,a_3,a_2,a_1)) \\
& + \nu(a_6,a_5,a_4,a_3,\mu^2(a_2,a_1)) \\
& - (-1)^{|a_1|} \nu(a_6,a_5,a_4,\mu^2(a_3,a_2),a_1) \\
& + (-1)^{(|a_2|+|a_1|)}\nu(a_6,a_5,\mu^2(a_4,a_3),a_2,a_1) \\
& - (-1)^{(|a_3|+|a_2|+|a_1|)} \nu(a_6,\mu^2(a_5,a_4),a_3,a_2,a_1) \\
& + (-1)^{(|a_4|+|a_3|+|a_2|+|a_1|)}\nu( \mu^2(a_6,a_5),a_4,a_3,a_2,a_1) .
\end{align*}
We check that 
\begin{align*} 
144 \mu^6(u,v,f_1,u,e_1, v) & = -9 f_0; \\
144\mu^6(f_1,u,v,u,e_1, v) & = 5f_0;\\
144\mu^6(f_1,u,e_1,v,u, v) & = 9 f_0;\\
144\mu^6(f_1,f_1,u,e_1,v, f_1) & = 11f_1.
\end{align*}
Hence
\begin{align*}  
9f_0 & = - \mu^2( \nu (u,v,f_1,u,e_1), v) - \mu^2( u, \nu (v,f_1,u,e_1,v)) - \nu(f_1,f_1,u,e_1,v) \\
& = - \nu(f_1,f_1,u,e_1,v) .
\end{align*}
Similarly, 
\begin{align*}
- 5f_0 & =  -\mu^2( \nu(f_1,u,v,u,e_1), v) - \mu^2( f_1, \nu(u,v,u,e_1,v)) + \nu(f_1,u,e_1,e_1,v) - \nu(f_1, f_1 ,u,e_1,v) \\
& =  \nu(f_1,u,e_1,e_1,v) + 9 f_0,
\end{align*}
so $\nu(f_1,u,e_1,e_1,v) = -14 f_0 $. Next,
\begin{align*}
-9f_0 &= - \mu^2( \nu (f_1,u,e_1,v,u), v) - \mu^2( f_1, \nu(u,e_1,v,u,v)) -
\nu(f_1,u,e_1,v,f_1) + \nu(f_1, u,e_1,e_1,v)  \\
&= -\nu(f_1,u,e_1,v,f_1) + \nu(f_1, u,e_1,e_1,v) \\
& =-\nu(f_1,u,e_1,v,f_1) - 14 f_0  
\end{align*}
so $\nu(f_1,u,e_1,v,f_1)=5f_0$. The \emph{coup de gr\^{a}ce} is the self-contradictory equation
\begin{align*}
- 11f_1 & =  -\mu^2( \nu(f_1,f_1,u,e_1,v), f_1) - \mu^2( f_1, \nu(f_1,u,e_1,v,f_1)) \\
& = 9 \mu^2(f_0,f_1) - 5 \mu^2(f_1,f_0) = -14 f_1 .
\end{align*}
This proves that $m^6\neq 0$.

\end{document}